\documentclass[ms,copyedit]{informs3TD} 

\OneAndAHalfSpacedXI 
\usepackage{caption}
\usepackage{subcaption}
\usepackage{algorithm} 
\usepackage{algorithmic} 
\newcommand{\bA}{ \mathbf{A} }

\newcommand{\bb}{ \mathbf{b} }

\newcommand{\bC}{ \mathbf{C} }
\newcommand{\bc}{ \mathbf{c} }
\newcommand{\bD}{ \mathbf{D} }

\newcommand{\bE}{ \mathbf{E} }

\newcommand{\bP}{ \mathbf{P} }

\newcommand{\bv}{ \mathbf{v} }

\newcommand{\bX}{ \mathbf{X} }
\newcommand{\bx}{ \mathbf{x} }
\newcommand{\bY}{ \mathbf{Y} }
\newcommand{\by}{ \mathbf{y} }

\newcommand{\FO}{ \text{FO} }
\newcommand{\IO}{ \text{IO} }
\newcommand{\bepsilon}{ \mathbf{\epsilon} }
\newcommand{\bz}{ \mathbf{z} }

\newcommand{\bzero}{\mathbf{0}}

\newcommand{\MLIO}{\text{MLIO}}
\newcommand{\MIBP}{\text{MLIO\textrm{-}MIBP}}

\newcommand{\SEQMLIO}{\text{SEQ\textrm{-}MLIO}}
\newcommand{\EMBMLIO}{\text{EMB\textrm{-}MLIO}}

\newcommand{\cD}{\mathcal{D}}

\usepackage{tikz}
\usetikzlibrary{arrows,calc,positioning}
\tikzstyle{intt}=[draw,text centered,minimum size=6em,text width=5.25cm,text height=0.34cm]
\tikzstyle{intl}=[draw,text centered,minimum size=2em,text width=2.75cm,text height=0.34cm]
\tikzstyle{int}=[draw,minimum size=2.5em,text centered,text width=3.5cm]
\tikzstyle{intg}=[draw,minimum size=3em,text centered,text width=6.cm]
\tikzstyle{sum}=[draw,shape=circle,inner sep=2pt,text centered,node distance=3.5cm]
\tikzstyle{summ}=[drawshape=circle,inner sep=4pt,text centered,node distance=3.cm]

\usepackage{natbib}
 \bibpunct[, ]{(}{)}{,}{a}{}{,}%
 \def\bibfont{\small}%
 \def\BIBand{and}%
 
\definecolor{chicago-maroon}{RGB}{128,0,0}

\usepackage{color}
\definecolor{hopkins-blue}{RGB}{0,45,114}
\definecolor{columbia-blue}{RGB}{185, 217, 235}
\definecolor{chicago-maroon}{RGB}{128,0,0}
\definecolor{cornell-red}{RGB}{179,27,27}
\definecolor{cmu-red}{cmyk}{0,1.00,0.79,0.20}
\definecolor{lawngreen}{RGB}{0,250,154}
\definecolor{gray}{RGB}{192,192,192}
\newcommand{\td}[1]{\textbf{\textcolor{cmu-red}{[Tinglong: #1]}}}

 \usepackage[colorlinks,citecolor=chicago-maroon,urlcolor=chicago-maroon,linkcolor=chicago-maroon]{hyperref}
\usepackage[nameinlink]{cleveref}

\usepackage{tikz}
\usepackage{subcaption}

\bibpunct[, ]{(}{)}{,}{a}{}{,}%
\def\bibfont{\small}%
\def\BIBand{and}%
\usepackage[colorlinks,citecolor=chicago-maroon,urlcolor=chicago-maroon,linkcolor=chicago-maroon]{hyperref}
\usepackage[nameinlink]{cleveref}
\crefname{assumption}{Assumption}{Assumptions}
\crefname{lemma}{Lemma}{Lemmas}
\crefname{theorem}{Theorem}{Theorems}
\crefname{corollary}{Corollary}{Corollaries}
\crefname{proposition}{Proposition}{Propositions}
\crefname{claim}{Claim}{Claims}
\crefname{definition}{Definition}{Definitions}
\crefname{subclaim}{Subclaim}{Subclaims}
\crefname{procedure}{Procedure}{Procedures}
\crefname{algorithm}{Algorithm}{Algorithms}
\crefname{example}{Example}{Examples}
\crefname{figure}{Figure}{Figures}
\crefname{section}{Section}{Sections}
\crefname{appendix}{Appendix}{Appendices}
\crefname{table}{Table}{Tables}
\crefname{equation}{}{}

\TheoremsNumberedThrough     
\ECRepeatTheorems

\EquationsNumberedThrough    

\MANUSCRIPTNO{MS-0001-1922.65}

\begin{document}



\RUNTITLE{A Preference-Aware Inverse Optimization Approach}
\RUNAUTHOR{Ahmadi, Dai, and Ghobadi}
\TITLE{You Are What You Eat: A Preference-Aware \\ Inverse Optimization Approach
}

\ARTICLEAUTHORS{%
\AUTHOR{Farzin Ahmadi{$^\ast$} \hspace{0.1in}  Tinglong Dai{$^\dagger$} \hspace{0.1in} Kimia Ghobadi{$^\ast$}
\medskip
\AFF{ 
{$^\ast$}Department of Civil and Systems Engineering and  Malone Center for Engineering in Healthcare, \\ Johns Hopkins University, Baltimore, Maryland 21218, \EMAIL{\{fahmadi1,kimia\}@jhu.edu}\\
\smallskip
{$^\dagger$}Carey Business School, Johns Hopkins University, Baltimore, Maryland 21202, \EMAIL{dai@jhu.edu}}
} 
}
\ABSTRACT{A key challenge in the emerging field of precision nutrition entails providing diet recommendations that reflect both the (often unknown) dietary preferences of different patient groups  and known dietary constraints specified by human experts. Motivated by this challenge, we develop a preference-aware constrained-inference approach in which the objective function of an optimization problem is not pre-specified and can differ across various segments. Among existing methods, clustering models from machine learning are not naturally suited for recovering the constrained optimization problems, whereas constrained inference models such as inverse optimization do not explicitly address non-homogeneity in given datasets. By harnessing the strengths of both clustering and inverse optimization techniques, we develop a novel approach that recovers the utility functions of a constrained optimization process across clusters while providing optimal diet recommendations as cluster representatives. Using a dataset of patients' daily food intakes, we show how our approach generalizes stand-alone clustering and inverse optimization approaches in terms of adherence to dietary guidelines and partitioning observations, respectively. The approach makes diet recommendations by incorporating both patient preferences and expert recommendations for healthier diets, leading to \emph{structural} improvements in  both patient partitioning and nutritional recommendations for each cluster.
An appealing feature of our method is its ability to consider \emph{infeasible} but informative observations for a given set of dietary constraints. 
The resulting recommendations correspond to a broader range of dietary options, even when they limit unhealthy choices.
}




\maketitle
\section{Introduction}
Precision nutrition, an emerging field that addresses ``the practical question of what to eat to stay healthy,''   has garnered significant attention in recent years, as imbalanced nutrition has emerged as a key contributor to numerous health issues that cost ``hundreds of billions of dollars'' each year \citep[][p. 735]{Rodgers2020}. In relation to this question, adhering to healthy dietary habits plays an essential role.  In practice, physicians often recommend specific diet regimens, for instance, the Dietary Approaches to Stop Hypertension (DASH) diet to control patients' sodium intake \citep{dash_diet_2020}. However, 
these diet regimens are rarely tailored to patients with their lifestyles or dietary preferences in mind, making long-term patient adherence challenging \citep{bazrafkan2021overweight, downer2016predictors,inelmen2005predictors}. Incorporating dietary preferences helps generate more palatable diet recommendations (i.e., preference-aware diet recommendations), which is a crucial building block of precision nutrition \citep{Rodgers2020}. Yet, dietary preferences, which are commonly reflected in each segment's objective functions, are rarely known; hence, broad and general dietary guidelines continue to be the norm.

\begin{figure}[!ht]
\begin{center}
\includegraphics[width =0.65 \linewidth]{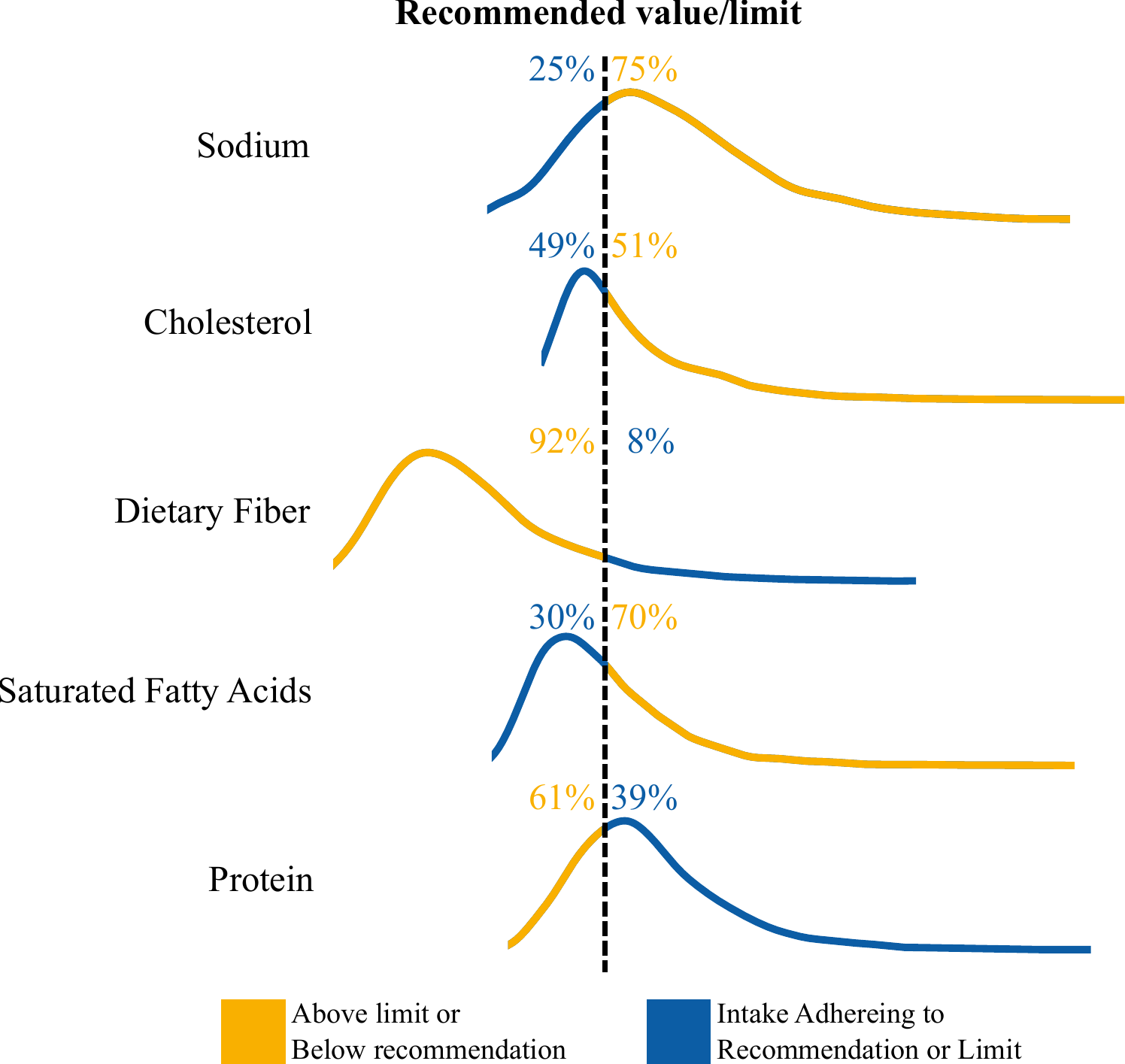}
\caption{\footnotesize A schematic figure representing the distribution of the nutritional intakes of the respondents from the NHANES dataset compared with the recommended values or limits set by the DASH dietary behavior plan. The majority of the respondents (at least half) exhibit unhealthy behaviors with regard to many essential nutrients, including sodium, cholesterol, dietary fiber, and saturated fatty acids.
} \label{Fig:Intro_sodium_figure}
\end{center}
\end{figure}

Interestingly, diet planning was among the pioneering application areas of operations research techniques, particularly linear programming \citep{dantzig1965linear, stigler1945cost}. In recent years, the focus of this area has shifted to using observed food-intake \emph{data} to help inform diet recommendations. However, commonly used data-driven approaches are more suitable for replicating patient behaviors (both desirable and non-desirable ones), which are often at odds with expert recommendations. As an illustrative example, in \cref{Fig:Intro_sodium_figure}, we plot the dietary behaviors of the National Health and Nutrition Examination Survey (NHANES) respondents, which includes nearly 10,000 respondents who record their daily food choices for two days.\footnote{\url{https://www.cdc.gov/nchs/nhanes/index.htm}} As the figure shows, the majority of the respondents do not adhere to the limits set forth by the DASH dietary guidelines for a majority of categories.
Clearly, a na\"ive clustering model can replicate the unhealthy behaviors that a significant portion of the patients exhibit. The presence of  diverse patient objectives (e.g., regulating calories, maximizing taste, and minimizing cost) and multiple dietary constraints (e.g., fiber intake should be within a certain range) further complicates the diet recommendation problem. In such circumstances, human intervention is essential to guide diet recommendation algorithms.

Existing machine learning algorithms are not immediately suitable for this type of problem. As powerful as these algorithms are, when used to ``recover'' optimization problems, they mostly operate under the assumption that the optimization problems are not subject to hard constraints and handle occasional constraints by introducing large penalties.
\footnote{One approach to address this major limitation of machine learning algorithms is ``constraint learning,'' that is, learning constraints from examples. Yet, constraint learning is a nascent field in machine learning, with limited evidence of its effectiveness and few real-world applications \citep{de2018learning}.}
Parallel to machine learning, inverse optimization focuses on recovering an optimization problem from a given set of observations \citep[see, e.g.,][]{Ahuja2001,aswani2018inverse,Chan2019,ghobadi2018robust}; theoretically, the original optimization (a.k.a. ``forward'' optimization) problem is recovered such that the observation becomes its optimal solution. 
However, a shortcoming of inverse optimization is that it usually applies to a single or a small number of observations \citep{Ahuja2001}. Although applying inverse optimization to the average or the worst case of multiple observations is possible \citep{esfahani2018data, ghobadi2018robust}, existing inverse optimization techniques  
assume all observations share the same feasible set and objective function.

\begin{figure}[t]
\begin{center}
\includegraphics[width =1 \linewidth]{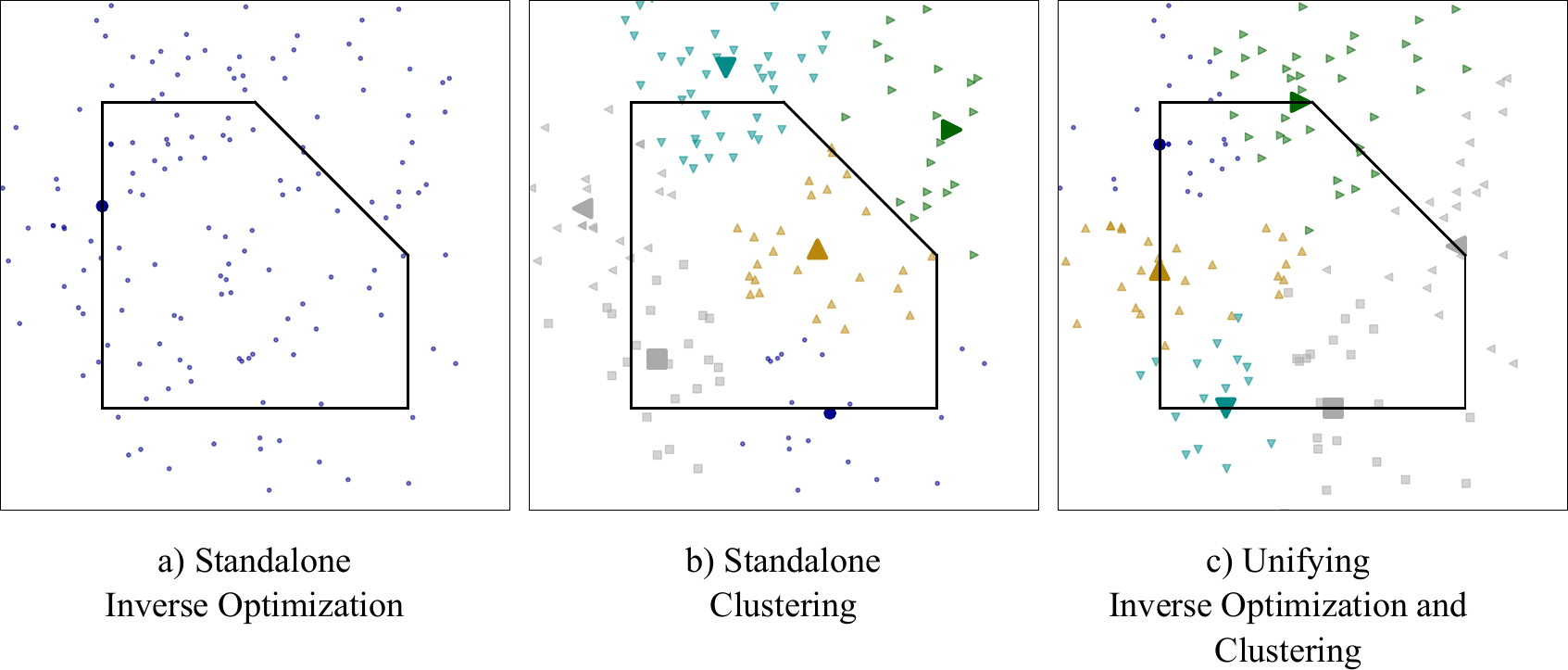}
\caption{\footnotesize 
A representative example of how inverse optimization combined with machine learning incorporates expert recommendation in clustering non-homogeneous observations, whereas both stand-alone inverse optimization and machine learning overlook important aspects, namely the differences between observed decisions and the known constraints of the problem, respectively. In the figure, smaller, similar shapes showcase decisions within the same cluster, whereas larger shapes show the cluster representative that the model learns. The shaded areas represent the partitioning.} \label{Fig:Model_descriptions_intro}
\end{center}
\end{figure}

\cref{Fig:Model_descriptions_intro} illustrates a bi-dimensional constrained environment where decisions are observed over potentially different objectives while a proposed shared feasible set is in place. As \cref{Fig:Model_descriptions_intro}(a) shows, when inverse optimization models are applied to recover objective vectors and optimal decisions, they generally assume all observations share the same feasible set and objective function. A stand-alone inverse optimization model is not developed with partitioning in mind; thus, in the face of all the observed decisions, stand-alone inverse optimization will yield a single optimal solution pertaining to a one-size-fits-all scenario. \cref{Fig:Model_descriptions_intro}(b), on the other hand, shows that stand-alone clustering  models do not consider the known constraints of the setting and result in suboptimal or even infeasible solutions. In other words, a na\"ive machine learning model is unable to capture the deviations and often provides recommendations that are---although similar to observed behaviors---undesirable (e.g., they tend to recommend diets high in sodium when patients' dietary choices tend to be excessively high in sodium). This tendency is in contrast to a combined model that incorporates the ability of machine learning models to capture non-homogeneity (i.e., potential variation across decision makers) and the ability of inverse optimization models to recover objective functions, which would provide better recommendations for each cluster of patients. 

Most notably, such a combined approach produces \emph{clusters} that differ from those produced by a stand-alone machine learning model in that each cluster contains a portion of the feasibility region with extreme points (which can be optimal for a realization of the unknown parameters of the optimization problem), whereas the stand-alone machine learning approach produces several interior clusters that are ``landlocked'' by other clusters. \cref{Fig:Model_descriptions_intro}(c) demonstrates how our approach yields \emph{structurally} different solutions in terms of both partitioning and the recommended solutions for each cluster.  Clearly, the combined approach is more than just layering an optimization problem on top of a clustering problem; rather, it allows for interaction between optimization and clustering.

In line with the spirit of \cref{Fig:Model_descriptions_intro}(c), we develop a novel methodology that unifies machine learning and inverse optimization (denoted as $\MLIO$) for constrained decision-making environments with unknown performance metrics and known constraints. Our methodology leverages the strengths of clustering models (in handling large and non-homogeneous datasets) and inverse optimization (in recovering the utility function under constraints and providing optimality guarantees) in large datasets.  
Using this methodology, we enforce desirable constraints on the cluster representatives and guarantee optimality conditions. Specifically, we formulate an optimization model that incorporates the known parameters of the constrained (forward) optimization problems and the set of all observed decisions. The model then optimizes over a loss criterion and recovers a number of clusters and the same number of (forward) optimization problems. We develop the model in such a way that for each of the recovered optimization problems and clusters, an optimal solution minimizing the loss towards the observed decisions within the cluster is also obtained. Furthermore, we present two solution approaches to the $\MLIO$ problem  for larger datasets. The first, \emph{sequential} approach treats the clustering and inverse optimization steps sequentially to provide a feasible solution for the general problem, whereas the second, \emph{embedded} approach builds on an initial partitioning by solving the inverse optimization model and reassigns observations jointly and iteratively. 

Additionally, we demonstrate that for the case of infeasible observed solutions, our proposed methodology has the desirable effect of providing optimal solutions that deviate incrementally from existing solutions, as shown in  \cref{Fig:Model_descriptions_intro}(c) for clusters including infeasible solutions. For our application, we specifically focus on the personalized diet recommendation problem. We show the performance and applicability of these solution approaches to the setting of the diet recommendation problem for a subset of patients from the NHANES dietary data.  We show  a na\"ive clustering model results in the replication of patients' existing dietary patterns that may reflect unhealthy lifestyles. By contrast, our MLIO model, through added flexibility of human inputs in the form of dietary constraints, recommends diets that adhere to DASH nutritional guidelines while incorporating dietary preferences inferred from the data. Even though they restrict unhealthy options, the resulting recommendations often correspond to a wider variety of dietary options.


The rest of this paper is organized as follows. In  \cref{sec:literature_review}, we review the related literature.  Next, we describe the problem setting from a theoretical standpoint in  \cref{sec:methodology} and detail our general $\MLIO$ model that unifies machine learning and inverse optimization. 
\cref{Sec:solution Approaches} discusses the solution method, where  \cref{sec:seq,sec:emb} detail and compare two efficient solution methods. Finally, \cref{sec:application} applies our novel approach to a diet recommendation problem and generates managerial implications. This paper concludes in \cref{sec:conclusion}.

\section{Literature} \label{sec:literature_review}
Our paper builds on and contributes to several streams of literature. 
Methodology-wise, we blend inverse optimization and unsupervised machine learning, so our paper is related to the literature on both approaches. Thematically, we contribute to the emerging field of precision nutrition and extend the scope of the healthcare operations management literature. More broadly, our paper is connected to the nascent literature on the human-algorithm connection. 

Inverse optimization is an inference tool for recovering optimization models  and has received increased attention since the seminal work by \cite{Ahuja2001}, who develop a method to recover the cost vector of a linear programming problem given the optimal solution of the problem.  Interest in inverse optimization first arose from finding parameters for combinatorial optimization problems such as the shortest-path problem  from given optimal solutions \citep{burton1992instance, zhang1999solution}.  In general, inverse optimization models aim to recover unknown parameters of an optimization problem by minimizing the loss in the optimality of the observed solutions. Early work in inverse optimization focuses on the case of a \emph{single} observed solution \citep[see, e.g.,][]{Ahuja2001,Chan2019, iyengar2005inverse, schaefer2009inverse}. To overcome this limitation, several recent papers \citep{aswani2018inverse,babier2021ensemble, ghobadi2020inferring, shahmoradi2021quantile} have extended this methodology to a limited number of (possibly noisy and/or suboptimal) observations. Recent inverse optimization works have shown the applicability of this approach in various health and non-health applications \citep{ahmadi2020inverse, Akhtar21,  aswani2019behavioral, barmann2018online, beil2003inverse,bertsimas2012inverse, chan2014generalized, chan2022inverse,chan2022inverse2,chow2012inverse, farago2003inverse, shahmoradi2022optimality, shahmoradi2021quantile}; we refer the reader to \cite{chan2021inverse} for a comprehensive survey of recent advances in theory and applications of inverse optimization. However, an underlying assumption across all these studies is that the observed decisions are solutions to  the same optimization problem. In other words,  all of the observed solutions are derived from the same utility function and the same feasible set. 
To the best of our knowledge, existing inverse optimization methods are not capable of distinguishing between observations that share the same feasible set but potentially different objectives and provide no means of partitioning given observations. This paper
bridges this gap in the literature. 

Various machine learning applications are also considered for constrained optimization problems. From this perspective, the majority of the literature focuses on methods to solve optimization problems (especially combinatorial optimization problems). The reader is referred to \citet{bengio2020machine} for a review of machine learning applications to continuous and combinatorial optimization problems. However, only a handful of studies have focused on learning solutions to optimization problems from given decisions. Among these studies, some map existing, noisy, and/or suboptimal observed decisions to better, near-optimal solutions based on the available knowledge of the problem context \citep{misra2018learning}. Deep neural network and reinforcement learning models have also been explored to generate ``good'' solutions for learning problems that involve decision making \citep{bengio2020machine,bastani2021learning}. However, most of the existing machine learning-based models do not guarantee the optimality of the mapped solutions and are prone to missing or violating important constraints in the optimization problem \citep{misra2018learning}, with a few exceptions: \citet{marquez2017imposing}, for example, explore the idea of imposing hard constraints on deep networks to guarantee optimality; as another example, \cite{elmachtoub2022smart} propose alternative loss functions in the prediction models that account for subsequent optimization problems. Machine learning techniques have also been used to infer the objective weights of optimization problems related to medical decision-making \citep{babier2018inverse,beam2018big}. 
Numerous studies investigate machine learning approaches and applications in healthcare services, for reviews of which we refer the reader to \cite{firdaus2018comparative} and  \cite{waring2020automated}.
Overall and mainly due to limitations from machine learning models to satisfy optimality of learned solutions in constrained data-driven settings, these models are not widely used for inference in constrained environments.

Taken together, neither machine learning (due to a general lack of optimality guarantees) nor inverse optimization (due to the inability to handle observations from different problems) alone is ideal for data-intensive constrained inference. 
However, by leveraging the strengths of both approaches, we develop a novel method that provides meaningful inference of unknown parameters and learns optimal solutions. The main contribution of our method is in its ability to recover different values for the unknown parameters of a general optimization problem through clustering a set of observations. We do so by using an augmented machine learning approach embedded with inverse optimization techniques that allow for optimally partitioning such observations. We show such an approach is superior to stand-alone applications of machine learning models in providing optimality guarantees for recommendations.

Our paper additionally contributes to the emerging stream of literature on human-algorithm connections.  \cite{dietvorst2015algorithm} define algorithm aversion as the tendency of forecasters and decision makers to use a human forecaster or algorithm against a statistical or evidence-based model, even when the evidence-based algorithms outperform human models. They also suggest that allowing users to modify algorithms mitigates algorithm aversion. This effect has been further examined in more recent research through incorporating considerations such as partial adherence to recommendations \citep{grand2022best}, treatment adherence \citep{lin2021does}, patients' resistance to medical artificial intelligence \citep{Longoni19}, physicians' reputation concerns \citep{Dai2020}, experts' updated belief in the algorithm's accuracy \citep{deVericourt2022}, and augmentation of algorithmic
decisions with human knowledge \citep{chen2022algorithmic}. Our paper is consistent with the literature in that we incorporate individuals' (unknown) dietary preferences when making diet recommendation decisions. Unlike previous research, which considers \emph{unconstrained} decision environments, our framework allows the algorithm's recommendations to explicitly reflect known decision constraints.

Our paper was motivated by a diet recommendation problem, which is among the earliest  application areas of operations research. The seminal work by \cite{stigler1945cost} models diet recommendation as an optimization problem and triggers broad interest among the optimization community \citep{dantzig1965linear}. The original diet recommendation problem entails finding the optimal intake amounts of different food items based on given constraints on food types and nutrients and given cost functions. Later studies have considered the diet recommendation problem at the individual \citep{maillot2010individual} and community levels \citep{buttriss2014diet,morgenstern2021perspective}. Other studies in the literature point out influential factors (e.g., convenience and taste) in dietary choices and note that recommendations should be focused on such aspects \citep{irz2016beyond}. Optimization models have been  used to model the diet recommendation problem as well \citep{gazan2018mathematical, ghobadi2018robust}. \cite{gazan2018mathematical} provide a review of diet optimization models that take both sustainability and acceptability into account. However, diet models usually suffer from recommendations that are not personalized for the population in question due to inadequate access to suitable cost functions or ill-defined feasible sets that are either infeasible or do not reflect the desirable features that the patients demand, or due to sole reliance on the observed behaviors in machine learning models that hinder adherence to healthier recommendations. Both inverse optimization models \citep{ahmadi2020inverse, ghobadi2018robust} and machine learning approaches \citep{ivancic2020diet} have considered the problem of providing meaningful and acceptable diets based on the observed decisions of patients. Although recent works consider providing clustering results with optimality guarantees for the preferences \citep{shahmoradi2022optimality} and consider different objective vectors for each observation \citep{birgestochastic}, to the best of our knowledge, our paper is the first to develop a hybrid approach to jointly (1) clustering individuals based on their dietary behaviors and (2) proposing optimal solutions and recovering unknown optimization parameters for each group.

Our work also contributes to the vibrant healthcare operations management literature, in which more tailored diet recommendations can be viewed as a gateway to improved quality of care. The literature in healthcare operations management is broad and ever-evolving. We refer the reader to  \cite{dai2020om}, \cite{keskinocak2020review}, and \cite{terwiesch2020om} for reviews of the healthcare operation management literature. Evidence within the healthcare operations management literature shows human judgment and expert recommendation could improve prediction' accuracy \citep{dai2021artificial, ibrahim2021eliciting}, and non-personalized treatment schedules are not optimal for different groups of patients \citep{suen2022design}. 



\section{A Preference-Aware Inverse Optimization Approach} \label{sec:methodology}
We consider a constrained decision-making environment where each decision maker (e.g., patient) is subject to a fixed and known feasible set. Different decision makers may have diverse objective functions even with shared decision constraints. Given a collection of observations, we develop a method that optimally partitions the observed decisions and recovers cost vectors in such a way that the representative decisions in each cluster are optimal for the recovered cost vector. This problem has parallels in applied settings such as the daily dietary behaviors of individuals who follow a particular diet (e.g., the DASH diet) and existing radiation therapy treatment plans for different patients~\citep{goldenberg2019new}. An appealing feature of our method is its ability to consider \emph{infeasible} observations for the given feasible set. Such observations may contain important information as the users may not be able or choose not to fully satisfy all constraints. For instance,  daily dietary behavior contains information about the user's preferences and palate, regardless of its feasibility for a particular diet.

We consider a set of different observations over a fixed and known polyhedral feasible set and assume the utility functions of the decision makers for the decisions are linear but unknown. The linearity assumption allows for a more tractable formulation while also allowing the objective function to be monotonous in terms of dietary preferences. Our approach consists of partitioning and learning: the partitioning component generates clusters of homogeneous decisions, and the learning component learns cluster representatives that optimally recover the unknown objective function based on observations within the cluster. 

In the rest of this section, we first lay the groundwork of our approaches in \cref{sec:methods_prelim} by setting up the modeling environment and specifying the learning component.  \cref{sec:methods_MLIO} presents the general problem of combing machine learning and inverse optimization. \cref{sec:bilinear_model} formulates the general problem as a mixed-integer bilinear program. 

\subsection{Preliminary: Inverse Optimization} \label{sec:methods_prelim}
We consider a setting where a set of potentially non-homogeneous observations is given along with a feasible set shared by all the observations. Corresponding to this setting is a forward linear optimization model $\FO$, where $\bc\in \mathbb{R}^n$ is the cost vector forming the utility function and $\Omega \subseteq \mathbb{R}^n$ is the feasible set where $\Omega \neq \emptyset$: 
\begin{subequations} \label{FO_i}
\begin{align}
\FO(\bc, \Omega): 
\underset{\bx}{\text{maximize}} & \quad  \bc' \bx \label{eq:IPP_obj}\\
\text{subject to} 
& \quad \bA \bx \geq \bb,  \label{eq:IPP_desConst}\\
& \quad \bx\in \mathbb{R}^n.
\end{align}
\end{subequations}
In the above formulation, $\bx \in \mathbb{R}^n$ represents the decision variables of $\FO$. The matrix $\bA \in \mathbb{R}^{m\times n}$ and the vector $\bb \in \mathbb{R}^{m}$ ($m$ is the number of constraints) constitute the required parameters to define the fixed feasible set $\Omega \subseteq \mathbb{R}^n$ of the forward optimization problem, and we have $\Omega = \left \{ x \in \mathbb{R}^n \mid \bA \bx \geq \bb \right \}$. Note the optimal solution set $\bX^{\ast} \subseteq \Omega$ for $\FO$ depends on the value of $\bc$. For a polyhedral feasible set $\Omega$, we denote the boundary of $\Omega$ as $\Omega^{opt}$, which contains all points in $\Omega$ that can be optimal for some cost vector $\bc \neq \bzero$. For a non-homogeneous set of observed decisions $\bX \subseteq \mathbb{R}^n$, our goal is to partition the given set of observed decisions into a given number of clusters $L$ such that the partitioning decision reflects the feasible set $\Omega$ and each cluster pertains to a group of decision-makers with the same objective function. Whereas clustering algorithms provide optimal partitioning (based on some metric) of the observations in the absence of constraints, we aim to incorporate the additional knowledge of the constraints into the partitioning scheme. To that end, we first discuss the inverse learning model as the specific inverse optimization approach we use to recover unknown cost vectors and learn optimal solutions to linear optimization problems. Then, we combine clustering and inverse optimization to recover optimization models for non-homogeneous observed decisions.

The learning component of our approach corresponds to an inverse optimization  model that is capable of recovering unknown parameters of optimization problems given homogeneous observed decisions. Because our goal is to generate clusters and cluster representatives that are optimal for the recovered utility functions, our inverse optimization model is the inverse learning framework, which is capable of learning optimal solutions. The core idea of inverse learning is to find a solution $\bz \in \mathbb{R}^n$ that results from minimal perturbation of given observed decisions $\bX_0$ attributed to $\FO$ to a singular optimal solution. By its nature of being an inverse optimization model, an inverse learning model is also capable of recovering a cost vector $\bc$ that makes the perturbed solution $\bz$ optimal for $\FO(\bc,\Omega)$. We denote by $\by \in \mathbb{R}^m$ the dual variables associated with constraints forming $\Omega$; $\bz \in \mathbb{R}^n$ the perturbed solution contained on $\Omega^{opt}$; $\bE \in \mathbb{R}^{n\times K}$ the perturbation matrix for the observed decisions; and $\epsilon^k$ be the $k^{th}$ column of $\bE$. Then, the inverse learning model is as follows: 
\begin{subequations} \label{IO_i}
\begin{align}
\IO({ \bX_0, \Omega}): 
\underset{\bc, \by, \bz, \bE}{\text{minimize}} & \quad \cD (\bz, \bX_0) \\
\text{subject to} 
& \quad \bA \bz \geq \bb,   \label{IOPrimalFeasblility1}\\ 
& \quad  \bc' \bz = \bb' \by ,   \label{SIOStrongDual}\\ 
& \quad \bA' \by  =  \bc, \label{IODualFeas1}\\ 
& \quad \bz = (\bx^k - \bepsilon^k ),   \quad \forall k \in \mathcal{K}\label{IOOnePoint}\\
& \quad \sum_{j =1}^m y_j =1,   \label{IORegularization}\\
& \quad \by \ge \bzero. \label{SIODualFeas2}
\end{align}
\end{subequations}

In the inverse learning model, the objective is to minimize some metric $\cD$ between the learned solution $\bz$ and the given observed decisions $\bX_0$. The constraints of $\IO$ guarantee the feasibility \eqref{IOPrimalFeasblility1} and optimality \eqref{SIOStrongDual} of $\bz$ for $\FO(\bc,\Omega)$. Note that although $\IO$ is non-convex, a growing body of literature proposes methods to solve $\IO$ by re-formulating it as a series of convex models \citep{aswani2018inverse, esfahani2018data, ghobadi2020inferring}. Additionally, it has been shown that  $\IO$ can be remodeled as a linearly constrained and convex optimization problem under mild assumptions on the constraint matrix $\bA$ \citep{ahmadi2020inverse,Chan2019}. \cite{ahmadi2020inverse} show that when $\Omega \neq \emptyset$, $\IO$ is feasible for any $\bX_0$ and for any metric $\cD$. We use this result in later sections to formulate a general model for non-homogeneous decisions $\bX$ that is always feasible when $\Omega \neq \emptyset$. 

Although $\IO$ is a powerful tool for recovering utility functions from observed decisions, it has an underlying assumption that $\bc$ is the same as all observed decisions and that $\bX_0$ is homogeneous. As such, by definition, $\IO$ is not capable of handling non-homogeneity in the observed decisions. In what follows, we consider clustering methodologies and develop a combined approach for learning a series of cost vectors and partitioning observed decisions.

\subsection{Recovering Optimization Models from Non-homogeneous Datasets} \label{sec:methods_MLIO}
In the presence of a non-homogeneous set of observed decisions, $\bX$, over multiple unknown utility functions $\bc_1,\hdots,\bc_L$, the practicality of the $\IO$ model discussed in the previous section is limited because $\IO$ is capable of learning a unified cost vector for all observed decisions. The presence of non-homogeneous decisions is analogous to a set of observed decisions from multiple decision makers/users over the same feasible set (e.g., the diet recommendation problem for different patient groups over the same dietary requirements). We aim to provide a method that can optimally (by minimizing a metric) partition the set of observed decisions by identifying the decisions that share the same decision maker and proposing a utility function for that decision maker. As such, we formulate a problem that simultaneously partitions the observed decisions into a given number of groups and recovers optimization models and their optimal solutions for each group.

Let $\bX \subseteq \mathbb{R}^{n}$ be a set of potentially non-homogeneous decisions made over $\Omega$. We aim to partition $\bX$ into a given number of clusters $\bX_1, \hdots , \bX_L$ while the model recovers utility parameters $\bc_1, \hdots , \bc_L$ and generates optimal solutions $\bx_1, \hdots , \bx_L$ for the recovered optimization problems $\FO_1(\bc_1,\Omega), \hdots , \FO_L(\bc_L,\Omega)$ such that $\bx_l$ are the optimal solutions contained in $\Omega^{opt}(\bc_l)\subseteq \mathbb{R}^n$, the set of all optimal solutions to $\FO(\bc_l,\Omega)$ $\forall l \in \mathcal{L} =  \left \{ 1, \hdots, L  \right \}$. We formulate this partitioning problem as follows: 
\begin{subequations} \label{MLIO}
\begin{align}
\MLIO({\bX, \Omega, L}): 
\underset{\bX_1,  \hdots , \bX_L, \bx_1,  \hdots  , \bx_L, \bc_1,  \hdots  , \bc_L }{\text{minimize}} & \quad \sum_{l=1}^{L} \hspace{0.05in} \mathcal{D} (\bx_l, \bX_l) \label{eq:opt_est_obj}\\
\text{subject to} 
& \quad \bx_l \in \Omega^{opt}(\bc_l), \quad \forall l \in \mathcal{L}  \label{eq:opt_est_const}\\
& \quad \bX_{\lambda} \cap \bX_{\gamma}  = \emptyset, \quad \forall \lambda,\gamma \in \mathcal{L}, \lambda \neq \gamma    \label{eq:partition_no_intersection} \\
& \quad \bigcup_{l = 1}^{L} \bX_l  = \bX,  \label{eq:partition} \\
& \quad \bc_l = \bzero \quad \text{iff} \ \Omega = \mathbb{R}^n. \quad \forall l \in \mathcal{L}  \label{eq:ML_well_defined} 
\end{align}
\end{subequations}
We refer to formulation \cref{MLIO} as the machine learning and inverse optimization ($\MLIO$) model. The MLIO model has four types of constraints: First, constraint \eqref{eq:opt_est_const} ensures learned solutions $\bx_l$ are optimal for $\FO(\bc_l,\Omega)$. Next, constraints \cref{eq:partition_no_intersection,eq:partition} impose the necessary criteria for $\left \{ \bX_1, \hdots, \bX_L  \right \}$ to be a partition of $\bX$.  In particular, constraint \eqref{eq:partition_no_intersection} ensures each observed decision will be included in exactly one cluster. Lastly, constraint \eqref{eq:ML_well_defined} ensures $\MLIO$ generalizes the unconstrained clustering problem. In the case of $\Omega = \mathbb{R}^n$,  by setting all $\bc_l = \bzero$, $\MLIO$ reduces to the unconstrained clustering problem, and the centroids that minimize $\cD$ over $\mathbb{R}^n$ are achieved. 

The objective of $\MLIO$ is to minimize the difference between the observed decisions in each cluster and the corresponding learned optimal solution for that cluster. In other words, MLIO looks for a partitioning of observed decisions that minimizes the total loss (as measured by the metric $\mathcal{D}$) for all clusters. We keep the metric $\cD$ generic and allow it to be specified by the human expert. To show $\MLIO$ is well defined, we first outline the feasibility of $\MLIO$ for different values of $L$.

\begin{proposition} \label{Prop:MLIO_feasibility}
$\MLIO({\bX, \Omega, L})$ is feasible for all $L \in \mathbb{N}$.
\end{proposition}

\cref{Prop:MLIO_feasibility} shows $\MLIO$ is well-defined for any given set of observed decisions; for $\Omega \neq \emptyset$, the feasible set of $\MLIO({\bX, \Omega, L})$, denoted by $\Psi_L$, is non-empty. An important characteristic of the $\MLIO$ model is that the feasible set of $\MLIO$ covers all possible partitions of $\bX$, denoted by $\bP_L(\bX)$. In other words, for any partition of $\bX$ to $L$ clusters, at least one solution for $\MLIO$ exists that results in that partition. This is formalized in \cref{Prop:MLIO_feasibility_all_partitions}.

\begin{proposition} \label{Prop:MLIO_feasibility_all_partitions} 
Let $\hat{\Psi}_L$ be the set of all partitions $(\bX_1,  \hdots , \bX_L) \in \bP_L(\bX)$ such that $\exists$ $ (\left\{ \bx_1,  \hdots  , \bx_L\right\},\left\{ \bc_1,  \hdots  , \bc_L\right\})$, where $(\bX_1,  \hdots , \bX_L, \bx_1,  \hdots  , \bx_L, \bc_1,  \hdots  , \bc_L) \in \Psi_L$. Then, $\hat{\Psi}_L$ covers $\bP_L(\bX)$.
\end{proposition}

\cref{Prop:MLIO_feasibility,Prop:MLIO_feasibility_all_partitions} establish that solving $\MLIO$ results in the optimal partition of $\bX$ with the desirable feature of recovering the appropriate number of optimization problems $\FO_1,\hdots,\FO_L$ and their corresponding optimal solutions $\bx_1,  \hdots, \bx_L$ and cost vectors $\bc_1,  \hdots, \bc_L$. Furthermore, the definition of $\MLIO$ provides a generalized approach for bridging the unsupervised machine learning (clustering) problem and the multi-observation inverse optimization problem.  For each cluster $\bX_l$, $\MLIO$ recovers a cost vector $\bc_l$ and a solution $\bx_l$ that is rendered optimal for $\FO(\bc_l,\Omega)$, so setting $\Omega = \mathbb{R}^n$ reduces $\MLIO$ to the unsupervised learning problem where all recovered cost vectors are equal to the zero vector and the learned solutions are the ones that minimize the metric over $\mathbb{R}^n$. On the other hand, setting $L = 1$ reduces $\MLIO$ to the inverse learning problem, where all observed decisions are assumed to emerge from the same optimization problem and the same cost vector $\bc$. 

\cref{Table:ML_IO_MLIO_comparison} compares (1) unsupervised machine learning approaches ($K$-means), (2) inverse optimization models ($\IO$), and (3) the hybrid model ($\MLIO$).   As the table shows, the main advantage of considering $\MLIO$ over existing inverse optimization models is $\MLIO$'s capability to handle potentially non-homogeneous decisions. As such, $\MLIO$ provides a general clustering approach that is aware of the decisions that are feasible for the problem ($\Omega$) where optimal clusters are learned based on the given metric $\cD$. The resulting optimal solution set of $\MLIO$ is a partition of the original observed decisions into $L$ groups, each group representing decisions that the formulations recognize as being for the same optimization problem and sharing the same cost vector. $\MLIO$ additionally provides optimal solutions for each such group in the partition.

\begin{table}[] 
\small
\begin{center}
        \caption{Comparison between unsupervised machine learning models, inverse optimization models, and $\MLIO$}
    \label{Table:ML_IO_MLIO_comparison}
    
    \begin{tabular} 
    {>{\centering}p{0.30\textwidth}|>{\centering}p{0.2\textwidth}>{\centering}p{0.2\textwidth}>{\centering\arraybackslash}p{0.20\textwidth}}
                                                                & \multicolumn{3}{c}{\textbf{}} \\      
{Features}        & Unconstrained Clustering \\  (e.g. $K$-means) & $\IO$ \\ {(e.g. Inverse learning)} & $\MLIO$  \\
\hline \hline
Proposing representative decisions    & $\begin{matrix} \\ \checkmark \\  \end{matrix}$ &  $\begin{matrix} \\ \checkmark \\  \end{matrix}$ &  $\begin{matrix} \\ \checkmark \\  \end{matrix}$    \\
\hline
Clustering capabilities    & $\begin{matrix} \\ \checkmark \\  \end{matrix}$ &  $\begin{matrix} \\  \\  \end{matrix}$ &  $\begin{matrix} \\ \checkmark \\  \end{matrix}$    \\
\hline
Considering non-trivial feasible sets ($\neq \mathbb{R}^n$)    & $\begin{matrix} \\  \\  \end{matrix}$ &  $\begin{matrix} \\ \checkmark \\  \end{matrix}$ &  $\begin{matrix} \\ \checkmark \\  \end{matrix}$    \\
\hline
Utility Function recovery    & $\begin{matrix} \\  \\  \end{matrix}$ &  $\begin{matrix} \\ \checkmark \\  \end{matrix}$ &  $\begin{matrix} \\ \checkmark \\  \end{matrix}$    \\
\hline
Learning optimal decisions    & $\begin{matrix} \\  \\  \end{matrix}$ &  $\begin{matrix} \\ \checkmark \\  \end{matrix}$ &  $\begin{matrix} \\ \checkmark \\  \end{matrix}$    \\
    \hline  \hline
    \end{tabular}
\end{center}
\end{table}
\normalfont

The definition provided in formulation \cref{MLIO} for $\MLIO$ is general and can be applied to any setting where inverse optimization can be modeled to learn optimal solutions. However,  how it can be explicitly modeled for different classes of $\FO$ problems is not obvious. As such, in what follows, we provide an explicit formulation for solving $\MLIO$ as a mixed-integer bilinear optimization problem that finds an optimal solution to $\MLIO$. We provide this formulation by building on the modeling approaches and techniques used for formulating $\IO$.

\subsection{A Mixed-Integer Bilinear Formulation} \label{sec:bilinear_model}
We now discuss how to explicitly model $\MLIO$ as an optimization problem and how to solve $\MLIO$ in a linear setting. $\MLIO$ can be modeled as a mixed-integer bilinear problem ($\MIBP$) in general by including binary variables indicating the inclusion of each observed decision in a cluster. We denote by $\bx_l \in \mathbb{R}^{n}$ the learned solution for cluster $l \in \mathcal{L}$ and $\bv \in \mathbb{R}^{K\times L}$ a binary matrix indicating the inclusion of observed decisions in each cluster, where $v_{k,l} = 1$ if and only if observed decision $k$ is in cluster $l$. Note that a one-to-one correspondence exists between the matrix $\bv$ and the partition in each feasible solution to $\MLIO$. We then have the following formulation:
\begin{subequations} \label{MLIO_MIBP}
\begin{align}
\MIBP({\bX, \Omega, L}): 
\underset{\bx, \bv, \bY, \bC, \bE }{\text{minimize}} & \quad  \sum_{l=1}^{L} \sum_{k=1}^{K} \hspace{0.05in} \mathcal{D} (\bx_l , \bx^k) \cdot v_{k,l}
\notag
\\
\text{subject to} 
& \quad \bA \bx_l \geq \bb  ,   \quad \forall l \in \mathcal{L} \label{MLIOPrimalFeasblility1}\\ 
& \quad  \bc_l' \bx_l = \bb' \by_l  ,   \quad \forall l \in \mathcal{L} \label{MLIOStrongDual}\\ 
& \quad \bA' \by_l  =  \bc_l,   \quad \forall l \in \mathcal{L} \label{MLIODualFeas1}\\ 
& \quad \bx_l \cdot v_{k,l} = (\bx^k - \bepsilon^{k} ) \cdot v_{k,l},   \quad \forall k \in \mathcal{K}, l \in \mathcal{L} \label{MLIOOnePoint}\\
& \quad \sum_{j = 1}^m y^j_l =1, \quad \forall l \in \mathcal{L}  \label{MLIORegularization}\\
& \quad \sum_{l=1}^{L} v_{k,l} = 1, \quad \forall k \in \mathcal{K} \label{MLIOclustering_constraint} \\
& \quad \by_l \ge \bzero, \quad \forall l \in \mathcal{L} \label{MLIODualFeas2} \\
& \quad \bv_{k,l} \in \left \{ 0,1 \right \}, \quad \forall k \in \mathcal{K}, l \in \mathcal{L}. \label{MLIO_membership}
\end{align}
\end{subequations}
In the above formulation, constraints \eqref{MLIOPrimalFeasblility1}, \eqref{MLIOStrongDual}, \eqref{MLIODualFeas1}, \eqref{MLIORegularization}, and \eqref{MLIODualFeas2} are similar to their definitions in \eqref{IO_i}. Constraints \eqref{MLIOOnePoint} ensure observed decisions sharing the same cluster are perturbed to a single solution. The constraints in formulation \eqref{MLIO_MIBP} ensure these perturbed solutions are contained in $\Omega^{opt} (\bc_l)$ for each cluster and, hence, are optimal for their respective recovered cost vectors. Constraints \eqref{MLIOclustering_constraint} ensure each observed decision is included in exactly one cluster. Note that for the special case of $\Omega = \mathbb{R}^n$, without loss of generality, we can remove constraints \eqref{MLIOPrimalFeasblility1}, \eqref{MLIOStrongDual},  \eqref{MLIODualFeas1}, \eqref{MLIORegularization}, and \eqref{MLIODualFeas2} because $\bA$ and $\bb$ do not exist. In this case, $\MIBP$ also reduces to the general clustering problem for the metric $\cD$ and $L$ number of clusters. Finally, note $\MIBP$ is developed as a generalization to the $\IO$ formulation to incorporate non-homogeneity in the observed decisions. \cref{prop:MIBP_generalization_IL} formalizes this statement.


\begin{proposition} \label{prop:MIBP_generalization_IL}
Let $(\bx^*, \bv^*, \bY^*, \bC^*, \bE^*)$ be optimal for $\MIBP(\bX, \Omega, L)$. Then, for all $l \in \mathcal{L}$, $(\bC^*_l, \bY^*_l, \bx^*_l, \hat{\bE^*_l})$ is optimal for $\IO (\bX^*_l,\Omega)$, where $\bx^*_l,\bC^*_l, \bY^*_l$ are the $l^{th}$ rows of $\bx^*,\bY^*, \bC^*$ and $\hat{\bE^*_l}$ is the matrix of rows of matrix  $\bE^*_l$ for which $v_{k,l} = 1$.
\end{proposition}

\cref{prop:MIBP_generalization_IL} shows how $\MIBP$ generalizes $\IO$ for non-homogeneous decisions. Instead of learning only one solution and recovering only one cost vector, $\MIBP$ is capable of learning $L$ optimal solutions and recovering $L$ cost vectors for each cluster. We note that although $\MIBP$ is a mixed-integer bilinear program, it can be reformulated as a mixed-integer linear program \citep{gupte2013solving}. 

Next,  in \cref{Thm:MLIO_Equal}, we show the important result that $\MIBP$ is equivalent to $\MLIO$. 
\begin{theorem} \label{Thm:MLIO_Equal}
$(\bx^*, \bv^*, \bY^*, \bC^*, \bE^*)$ is optimal for $\MIBP(\bX_0, \Omega, L)$ if and only if $\exists$ $ \bar{\bX} = \left \{ \bX_1,\hdots,\bX_L \right \} \in \bP_L(\bX)$, such that $(\bX,\bx^*, \bC^*)$ is optimal for $\MLIO(\bX, \Omega, L)$.
\end{theorem}

The $\MIBP$ formulation detailed in this section provides an explicit formulation to solve $\MLIO$ and learn the optimal partition of observed decisions in a general, potentially constrained environment. However, recall from \cref{sec:methods_prelim} that the presence of a combination of binary variables and bilinear or non-convex constraints in $\MIBP$ makes this model computationally expensive in the general case. Thus, in the following section, we provide two computationally efficient solution approaches to $\MLIO$. 

\section{Solution Approaches: Sequential vs. Embedded} \label{Sec:solution Approaches}

In this section, we introduce two solution approaches to $\MLIO$. 
\cref{sec:seq} discusses a \emph{sequential} solution approach where we first cluster the observed decisions independently from the given feasible set (equivalent to solving the unsupervised learning problem) and use $\IO$ to find optimal solutions and recover optimization problems for each cluster. Then, \cref{sec:emb} discusses an iterative, \emph{embedded} solution approach that finds clusters by minimizing the loss of points that are optimal for $\Omega$ from each cluster. The partition that is found by minimizing a convergence gap is then returned by the algorithm. Finally, \cref{sec:Numerical_Example} compares the two solution approaches in terms of optimal objective and run time with $\MLIO$. 

\subsection{A Sequential Approach to Solving $\MLIO$} \label{sec:seq}
\begin{figure}[t]
\begin{center}
\includegraphics[width =1 \linewidth]{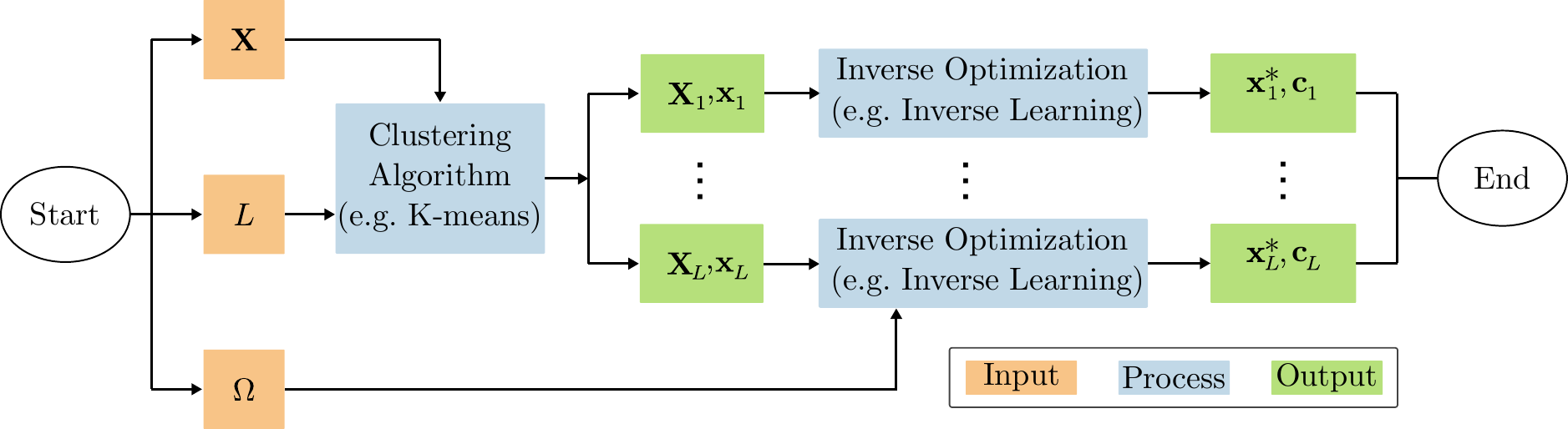}
\caption{\footnotesize Overview of the $\SEQMLIO$ methodology. This approach clusters the observations based on a given loss metric and finds inverse optimal solutions for each cluster.} \label{Fig:Diagram1}
\end{center}
\end{figure}

Our first solution approach to tackling the complexity of $\MIBP$ is to separate the clustering component from the learning component. In other words, we provide an approach that utilizes clustering models, such as the $K$-means method, to perform the clustering component of the learning process and use the $\IO$ model on each of the resulting clusters to learn optimal solutions. This approach provides an efficient method for finding a feasible solution to $\MLIO$ because inverse linear optimization models can be broken down into linearly constrained optimization problems. \cref{Fig:Diagram1} provides a schematic diagram of this approach. As shown in the figure, the non-homogeneous observed decisions $\bX$ undergo a clustering scheme using an unsupervised machine learning model, and the resulting $L$ clusters are each used as inputs to the $\IO$ model, alongside the known feasible set $\Omega$ one by one, resulting in recovering at most $L$ forward optimization problems and learning at most $L$ optimal solutions. We denote each learned optimal solution $\bx_l$ as a sequential machine learning and inverse optimization ($\SEQMLIO$) representative of cluster $l$. Whereas the $\SEQMLIO$ approach is computationally efficient, it bears the assumption that decisions are not influenced by the given feasible set and decisions closer to each other belong in the same cluster, which might not hold in highly constrained settings. Nevertheless, the $\SEQMLIO$ approach serves as a useful benchmark that can readily find a feasible solution for $\MLIO$. 
Proposition \ref{Prop:SEQ_MLIO_feasible_MLIO} states that the resulting solution of the approach is a feasible solution to $\MLIO$.

\begin{proposition} \label{Prop:SEQ_MLIO_feasible_MLIO}
The solution $(\left\{ \bX^1_0, \dots, \bX^L_0\right\}, \left\{\bx^*_1, \dots, \bx^*_L \right\}, \left\{\bc_1, \dots, \bc_L \right\})$ resulting from the $\SEQMLIO$ approach is a feasible solution for $\MLIO$.
\end{proposition}

\cref{Prop:SEQ_MLIO_feasible_MLIO} shows the solution obtained by
$\SEQMLIO$ is feasible for $\MLIO$. Yet, $\SEQMLIO$ does not provide a quality guarantee, because by design, it does not incorporate the available knowledge on $\Omega$ in the partitioning process, which plays a vital role in determining similarities among decisions.  
Furthermore, under $\SEQMLIO$, the resulting partitioning suffers from the same assumptions that rudimentary machine learning models are based on, which do not hold in general when $\Omega \neq \mathbb{R}^n$. Therefore, the $\SEQMLIO$ partition of $\bX_0$ is indifferent to $\Omega$, and additional available information on $\Omega$ does not affect the resulting partition and only affects the learned optimal solutions. To address this issue, in what follows, we consider another solution approach that captures $\Omega$ in the clustering algorithm.

\subsection{An Embedded Approach to Solving $\MLIO$} \label{sec:emb}

 To improve on the solution found by $\SEQMLIO$, one can either improve the partition or find cluster representatives that decrease the overall loss.  Because $\IO$ always learns cluster representatives with the minimum loss to the observations within a cluster, finding improved partitions is a natural direction to reduce the total loss. In this section,  we provide a modified heuristics method that is capable of providing partial optimality guarantees.
 
 To start, consider the feasible set of $\MLIO$ that contains all possible partitions of $\bX$ for a given number of clusters. Developing a heuristic that learns optimal or local optimal solutions is a daunting task. For such problems, previous literature has considered the notion of \emph{partial optimal} solutions as defined in \cref{Def:partial_optimal_solution}
 \citep{chakraborty2017k,selim1984k}.

\begin{definition} \label{Def:partial_optimal_solution}
Let $\bD (\bX_1,  \hdots , \bX_L, \bx_1,\hdots, \bx_L) = \sum_{i=1}^{L} \hspace{0.05in} \cD (\bx_i, \bX^i)$. 
The solution ($\bX_1,  \hdots , \bX_L, \bx_1,\hdots, \bx_L,\bc_1,\hdots, \bc_L $) is a partial optimal solution for $\MLIO$ if it satisfies the following:
\begin{enumerate}
    \item[(i)] $\bD (\bX_1,  \hdots , \bX_L, \bx_1,\hdots, \bx_L) \leq \bD (\bX'_1,  \hdots , \bX'_L, \bx_1,\hdots, \bx_L)$ for any  $\left \{ \bX'_1,\hdots,\bX'_L \right \} \in \bP_L(\bX)$.
    \item[(ii)] $\bD (\bX_1,  \hdots , \bX_L, \bx_1,\hdots, \bx_L) \leq \bD (\bX_1,  \hdots , \bX_L, \bx'_1,\hdots, \bx'_L)$ for any solution set $\left\{ \bx'_1,\hdots, \bx'_L \right\} \subseteq \Omega^{opt}$.
\end{enumerate}
\end{definition}

\cref{Def:partial_optimal_solution} provides a measure of quality for a proposed solution to $\MLIO$. A solution is partial optimal if the sub-problems of finding the best set of clusters with fixed cluster representatives (which corresponds to part (i) of the definition) and finding the best cluster representatives with a fixed set of clusters do not yield better solutions (which corresponds to part (ii) of the definition) in terms of the overall metric $\bD$.

Using \cref{Def:partial_optimal_solution}, we propose an algorithm that converges to a partial optimal solution for $\MLIO$. To this end, we embed the inverse learning model (which provides optimality guarantees for the learned parameters and solutions) into the clustering model (which handles the partitioning decisions) to return optimal solutions as cluster representatives. Such an approach will allow the model to update the partition of $\bX_0$ in a direction that reduces the total loss. Instead of finding point estimators in $\mathbb{R}^n$ that minimize loss to the observed decisions, our proposed embedded algorithm aims to find solutions contained in $\Omega^{opt}$ and reduce the distance between found solutions and the decisions in the clusters until a local optimum is reached. We show the algorithm concludes after a finite number of iterations and discuss desirable properties regarding the final output of the algorithm.  We then compare the results of these solution models in terms of optimal objective and run time with $\MLIO$ in \cref{sec:Numerical_Example,sec:application}.

\begin{figure}[t]
\begin{center}
\includegraphics[width =1 \linewidth]{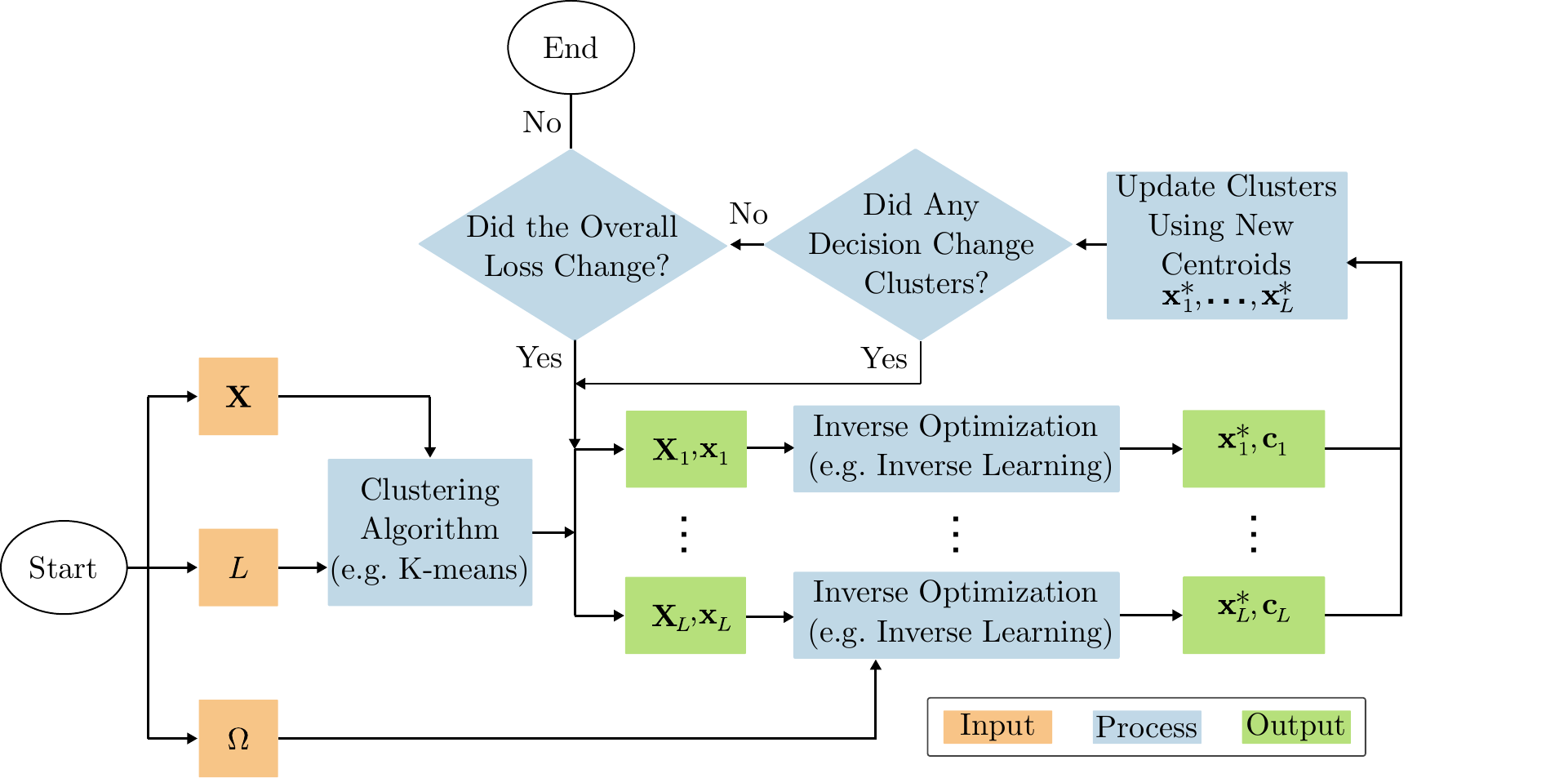}
\caption{\footnotesize Overview of the $\EMBMLIO$ methodology. This approach clusters the observations based on some loss metric and finds inverse optimal solutions for each cluster, whereas an embedded inverse optimization model guides the clustering updates in each iteration until the termination conditions are satisfied.} \label{Fig:Diagram2}
\end{center}
\end{figure}

The embedded approach $\EMBMLIO$ is depicted in \cref{Fig:Diagram2} similar to the $\SEQMLIO$ in \cref{sec:seq}. As indicated in \cref{Fig:Diagram2}, additional steps for the $\EMBMLIO$ approach can improve the quality of the final solution relative to the $\SEQMLIO$ approach. Details of the $\EMBMLIO$ approach are provided in \cref{e_algorithm_for_MLIO}. The algorithm starts with a distance-based clustering of the observations (e.g., similar to the $K$-means algorithm) and iteratively solves the $\IO$ problem for each cluster and updates the clusters based on the metric value of each observed decision to the newly learned optimal solutions. The algorithm terminates if no changes are made to the clusters  or if the total loss does not decrease. We first outline the feasibility of the learned solutions from \cref{e_algorithm_for_MLIO} for $\MLIO$ and show that for any observed decision set $\bX_0$, $\EMBMLIO$ returns a solution with a smaller total loss ($\cD$) relative to $\SEQMLIO$. 

\begin{proposition} \label{Prop:EMB_MLIO_feasible_MLIO}
The solutions $(\bX^1_0, \dots, \bX^L_0$, $\bx^*_1, \dots, \bx^*_L, \bc_1,\hdots, \bc_L)$ from each iteration of \cref{e_algorithm_for_MLIO} are feasible for $\MLIO$.
\end{proposition}

\cref{Prop:EMB_MLIO_feasible_MLIO} states that the solution returned using $\EMBMLIO$ (\cref{e_algorithm_for_MLIO}) contains optimal solutions for each cluster. In the following proposition, we go one step further by showing $\EMBMLIO$ always provides a solution that reduces the total loss in comparison to the solution obtained from $\SEQMLIO$.

\begin{proposition} \label{Prop:EMB_MLIO_vs_SEQ_MLIO}
Let $(\left\{ \bX^S_1, \dots, \bX^S_L\right\}, \left\{\bx^S_1, \dots, \bx^S_L \right\} , \left\{\bc^S_1, \dots, \bc^S_L \right\})$ be the $\MLIO$ feasible solution resulting from the $\SEQMLIO$ approach. Then, $ \exists$ $ (\left\{ \bX^E_1, \dots, \bX^E_L\right\}, \left\{\bx^E_1, \dots, \bx^E_L \right\}, \left\{\bc^E_1, \dots, \bc^E_L \right\})$ can be achieved from the $\EMBMLIO$ approach using \cref{e_algorithm_for_MLIO} such that 
$\sum_{l=1}^{L} \mathcal{D} (\bx^E_l, \bX^E_l) \leq \sum_{l=1}^{L} \mathcal{D} (\bx^S_l, \bX^S_l)$.
\end{proposition}

\begin{algorithm} 
	\caption{$\EMBMLIO$ for Clustering  Decisions}\label{e_algorithm_for_MLIO}
	\begin{algorithmic}[1]
	\STATE \textbf{INPUT:} $\bX$, $\Omega$, $L$
	\STATE \textbf{OUTPUT:} $\bX_1, ... , \bX_L$, $\bx_1, ... , \bx_L$, $\bc_1, ... , \bc_L$
	\STATE Select Initial Partition $\bX_1, ... , \bX_L$
	\STATE \textbf{Step:}
	\FOR{$i \in \left \{ 1,...,L \right \}$}
	    \STATE $\bc_{i}, \by_{i}, \bx_{i}, \bE_{i} \gets $Solve $\IO_i(\bX_i, \Omega)$
	\ENDFOR
	\FOR{ $\bx \in \bX$}
        \STATE Find $i^* = \arg\min \left \{ \cD(\bx, \bx_i) \quad \forall i \in \mathcal{L}\right \}$
        \STATE Add $\bx$ to $\bX^{New}$
    \ENDFOR
    \STATE $\bar{\bX} \gets \left \{ \bX_1,\hdots,\bX_L \right \}$, $\bar{\bX}^{New} \gets \left \{ \bX^{New}_1,\hdots,\bX^{New}_L \right \}$
    \IF {$\bar{\bX}^{New}$ = $\bar{\bX}$}
        \STATE $\bX_1, ... , \bX_L \gets \bX^{New}_1, ... , \bX^{New}_L$
        \STATE \textbf{STOP}
    \ENDIF
    \IF{$\bD (\bX^{New}_1,  \hdots , \bX^{New}_L, \bx_1,\hdots, \bx_L) = \bD (\bX_1,  \hdots , \bX_L, \bx_1,\hdots, \bx_L)$}
        \STATE $\bX_1, ... , \bX_L \gets \bX^{New}_1, ... , \bX^{New}_L$
        \STATE \textbf{STOP}
    \ELSE
        \STATE $\bX_1, ... , \bX_L \gets \bX^{New}_1, ... , \bX^{New}_L$
        \STATE Go to \textbf{Step}.
    \ENDIF

	\end{algorithmic}
\end{algorithm}

$\EMBMLIO$ is capable of finding local improvements, if any are available, when the solution from $\SEQMLIO$ is input to \cref{e_algorithm_for_MLIO}. Having established that \cref{e_algorithm_for_MLIO} is well-defined through \cref{Prop:EMB_MLIO_vs_SEQ_MLIO,Prop:EMB_MLIO_feasible_MLIO}, we now turn our attention to explaining the convergence and quality of the final solution of \cref{e_algorithm_for_MLIO}. We reiterate that because \cref{e_algorithm_for_MLIO} is a $K$-means type algorithm, one can prove it terminates in a finite number of iterations. To prove this result,  we first show in \cref{Lemma_EMBMLIO} that each possible partition of $\bX$ is visited at most once through the course of \cref{e_algorithm_for_MLIO}. 

\begin{lemma} \label{Lemma_EMBMLIO}
Any $ \bar{\bX} = \left \{ \bX_1,\hdots,\bX_L \right \} \in \bP_L(\bX)$ is visited at most once through in \cref{e_algorithm_for_MLIO}.
\end{lemma}

Using the results of \cref{Lemma_EMBMLIO} and the fact that \cref{e_algorithm_for_MLIO} is strictly decreasing in terms of total loss between iterations, we show in \cref{prop:ALg_1_conv_iterations} that \cref{e_algorithm_for_MLIO} converges to a solution in a finite number of iterations. We note the process of proving the convergence of \cref{e_algorithm_for_MLIO} is similar to the work of \cite{selim1984k}, who show similar results for the $K$-means algorithm. 

\begin{proposition} \label{prop:ALg_1_conv_iterations}
 \cref{e_algorithm_for_MLIO} converges in $\rho \leq \left| \bP_L(\bX)\right|$ iterations.
\end{proposition}

\cref{prop:ALg_1_conv_iterations} shows \cref{e_algorithm_for_MLIO} always converges to a solution in a finite number of iterations, which is capped at the number of possible partitions corresponding to the desirable number of clusters $L$. Note this result per se does not guarantee the performance of the final solution of \cref{e_algorithm_for_MLIO}. Next, in \cref{EMB_MLIO_convergence}, we show the output of the algorithm is a partial optimal solution for $\MLIO$.

\begin{theorem} \label{EMB_MLIO_convergence}
Let $(\bX^1, ... , \bX^L, \bx_1, ... , \bx_L, \bc_1, ... , \bc_L)$ be the output of \cref{e_algorithm_for_MLIO} on $(\bX_0,\Omega,L)$. Then, $(\bX^1, ... , \bX^L, \bx_1, ... , \bx_L, \bc_1, ... , \bc_L)$ is a partial optimal solution for $\MLIO$.
\end{theorem}

Based on the above discussions, the resulting partition from \cref{e_algorithm_for_MLIO} has the desirable characteristic of being a partial optimal solution for $\MLIO$. Additionally, using arguments similar to those in the literature \citep{selim1984k}, one can show the final $\MLIO$ solution found from \cref{e_algorithm_for_MLIO} is a local optimal solution of $\MLIO$ for a polyhedral $\Omega$. Stated differently, $\EMBMLIO$ is capable of efficiently providing an $\MLIO$ solution with desirable characteristics. 

\subsection{Numerical Illustration} \label{sec:Numerical_Example}
We now illustrate our proposed models through a bi-dimensional numerical example. We consider a fixed feasible set $\Omega$ representing known constraints for healthy lifestyles, and assume a number of observed decisions scattered over this fixed feasible set. We then compare how different models perform in partitioning these given decisions, learning optimal solutions, and recovering optimization models for these observed decisions. 

\begin{figure}[t]
\centering
\begin{subfigure}{0.3\textwidth}
    \includegraphics[width=\textwidth]{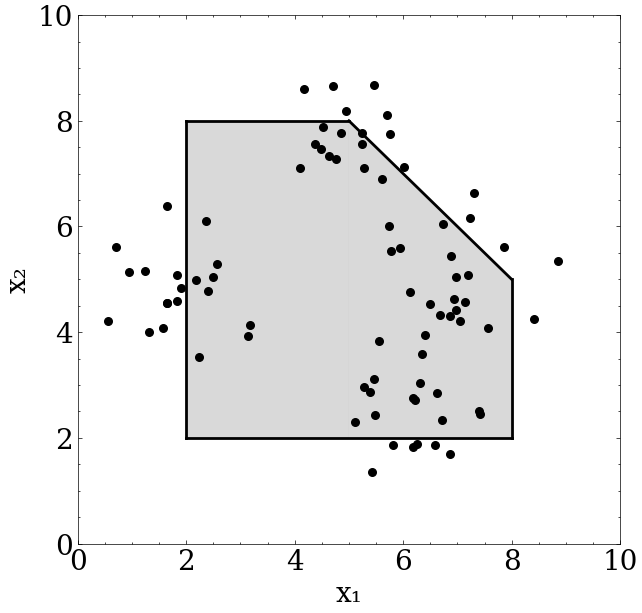}
    \caption{Feasible Set}
    \label{fig:2D_Partition_comparison_feasible_set}
\end{subfigure}
\hfill
\begin{subfigure}{0.3\textwidth}
    \includegraphics[width=\textwidth]{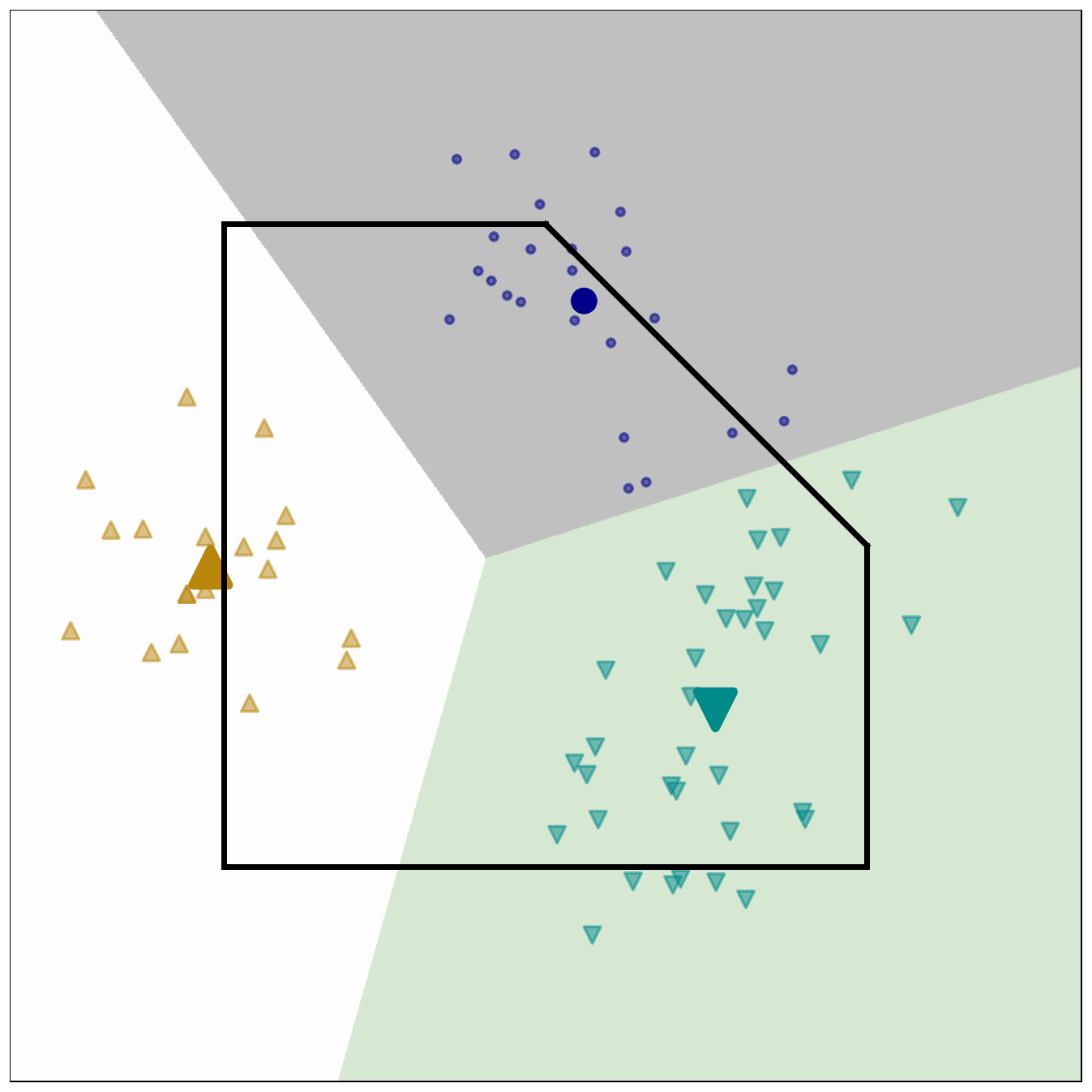}
    \caption{$K$-means}
    \label{fig:2D_Partition_comparison_kmeans}
\end{subfigure}
\hfill
\begin{subfigure}{0.3\textwidth}
    \includegraphics[width=\textwidth]{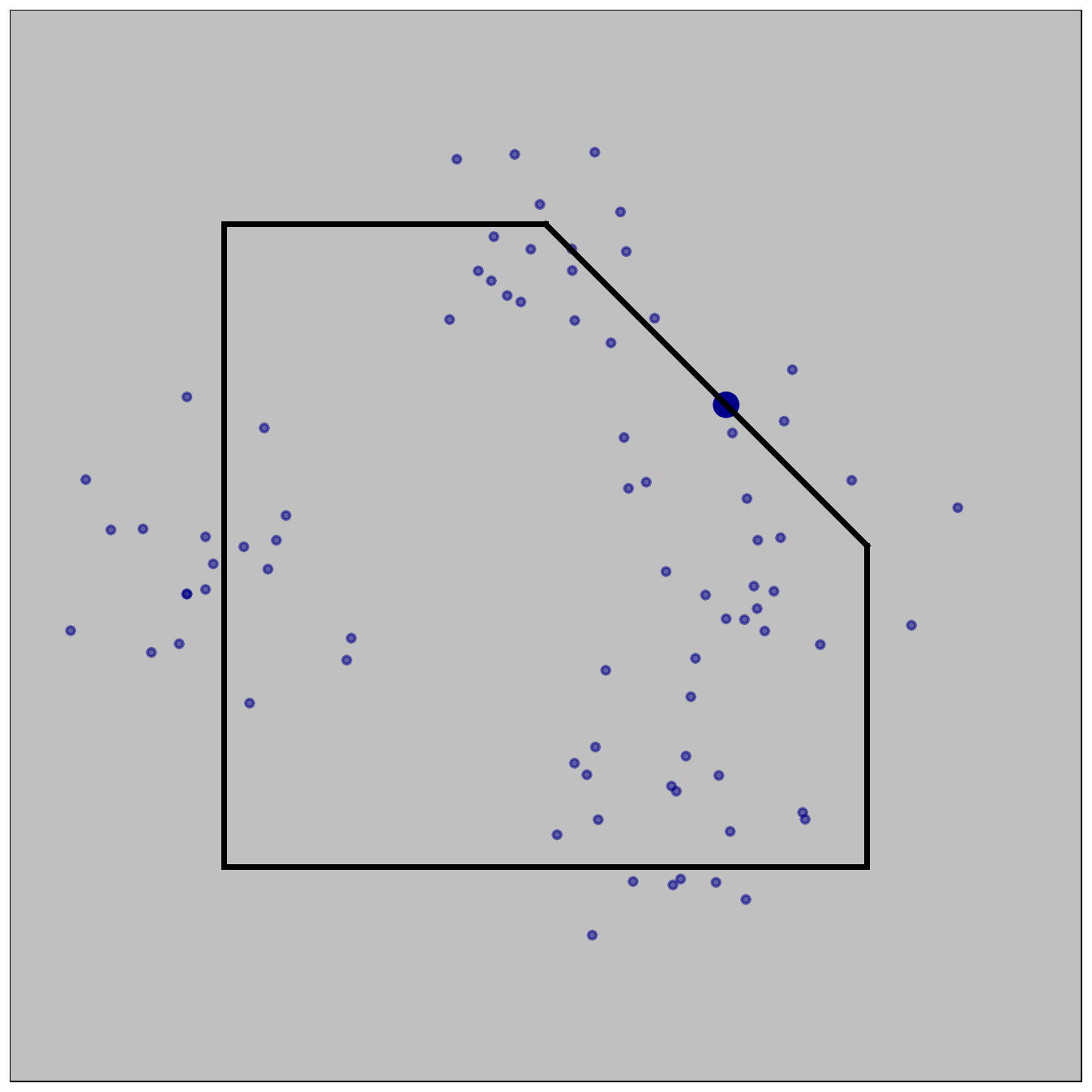}
    \caption{$\IO$}
    \label{fig:2D_Partition_comparison_IL}
\end{subfigure}
\hfill
\begin{subfigure}{0.3\textwidth}
    \includegraphics[width=\textwidth]{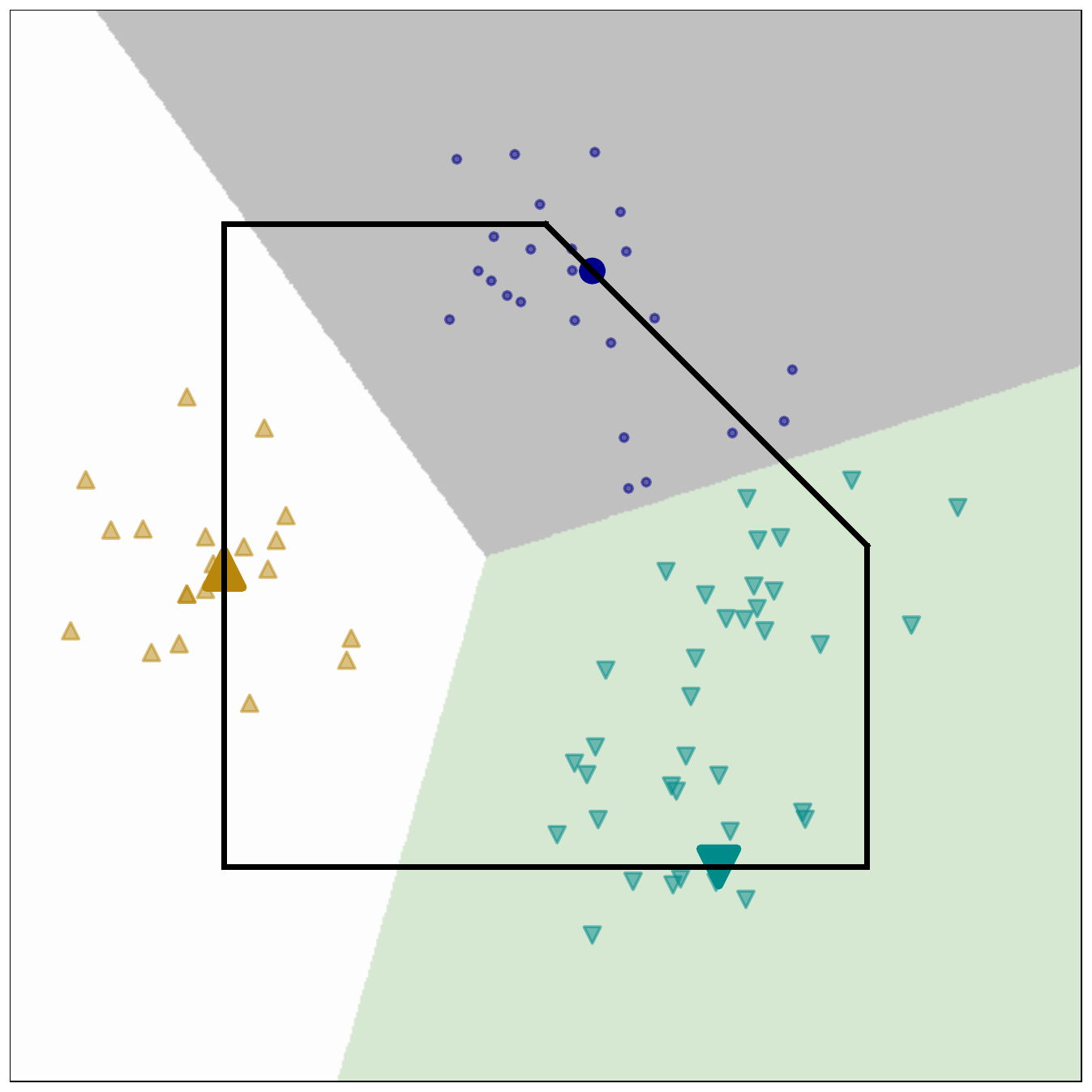}
    \caption{ $\SEQMLIO$}
    \label{fig:2D_Partition_comparison_SEQ}
\end{subfigure}
\hfill
\begin{subfigure}{0.3\textwidth}
    \includegraphics[width=\textwidth]{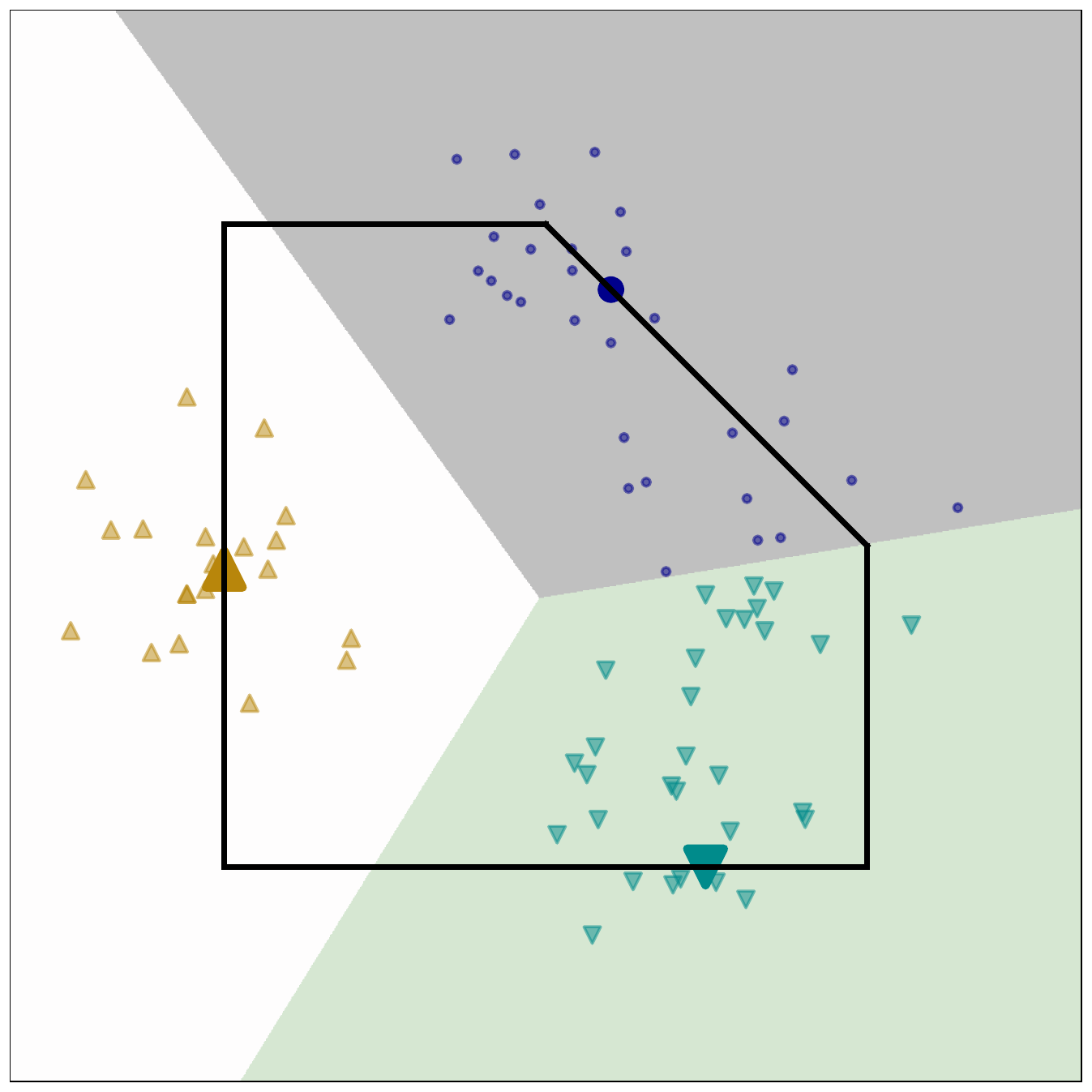}
    \caption{$\EMBMLIO$}
    \label{fig:2D_Partition_comparison_EMB}
\end{subfigure}
\hfill
\begin{subfigure}{0.3\textwidth}
    \includegraphics[width=\textwidth]{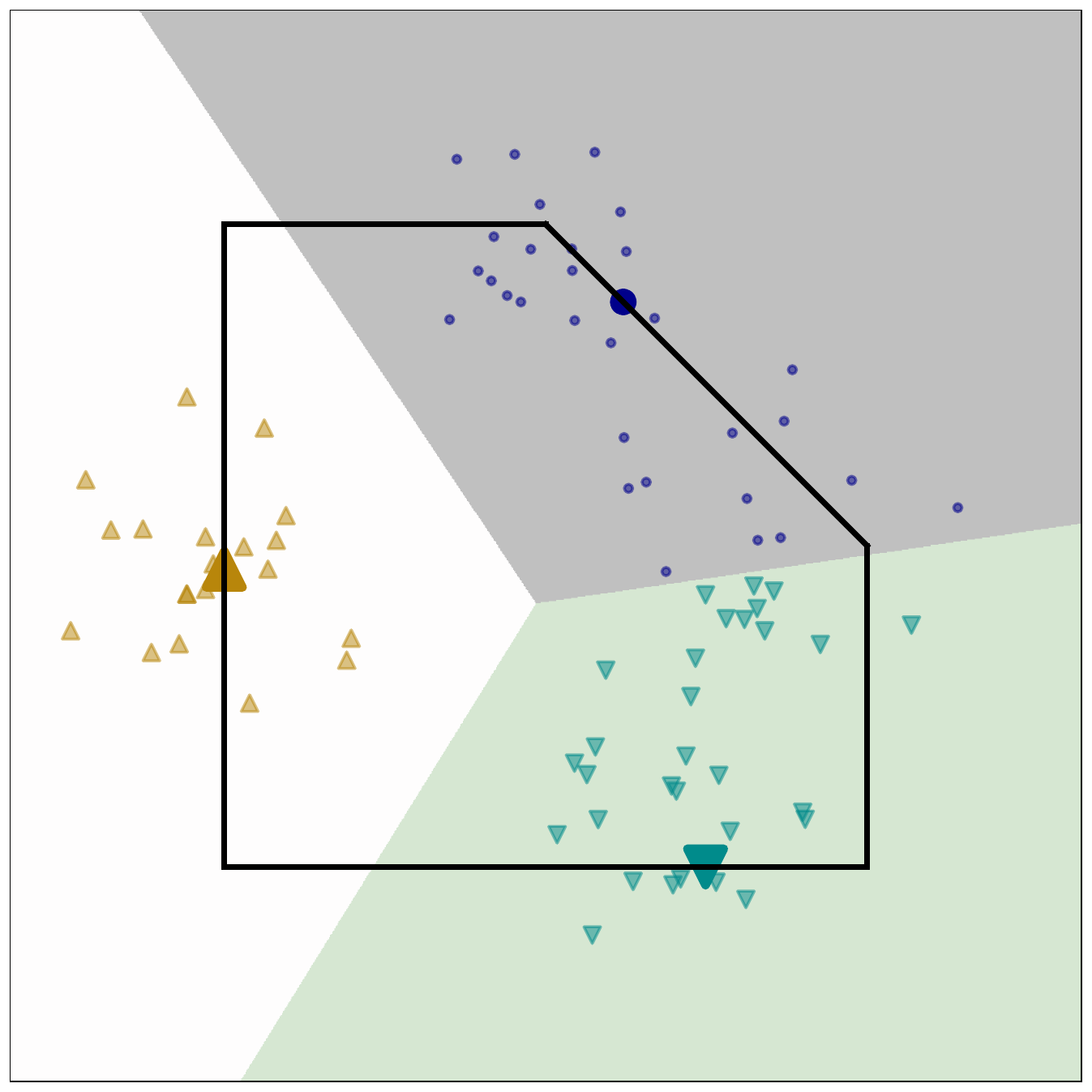}
    \caption{$\MIBP$}
    \label{fig:2D_Partition_comparison_MIQCP}
\end{subfigure} \vspace{0.2in}
\caption{Comparison of partitioning results using $\MIBP$, $\SEQMLIO$ and $\EMBMLIO$ as well as stand-alone inverse optimization and clustering models for a two-dimensional example problem. Combining machine learning and inverse optimization approaches allows $\MIBP$ and $\EMBMLIO$ to modify the color mapping of clusters and representative solutions to adhere to knowledge gained using $\Omega$.}
\label{fig:2D_Partition_comparison}
\end{figure}

\begin{example} \label{example:2DLinearOptimization_example1}
Let $\Omega=\left \{ (x_1,x_2)~|~  1 \leq x_1 \leq 5, ~~~ 1 \leq x_2 \leq 10,  ~~~ x_1 + x_2 \leq 15 \right \}$ be a bi-dimensional polyhedral feasible set of a constrained decision making environment with unknown cost vectors. This fixed and known feasible set represents the constraints for healthy lifestyles that $\MLIO$ models use to provide a partitioning of the observed decisions for a given number of clusters. Additionally, a set of 80 observed decisions was randomly generated that included both feasible and infeasible points for $\Omega$. The observed decisions and $\Omega$ are shown in \cref{fig:2D_Partition_comparison}(a).  For linear problems, the relationship between optimality and vicinity to the boundary of $\Omega$ will cause the $\MLIO$ models to provide an alternative partitioning in comparison to stand-alone unsupervised clustering models, as shown in \cref{fig:2D_Partition_comparison}.
\end{example}

\cref{fig:2D_Partition_comparison} and \cref{Table:2D_example_comparison} compare the results of employing (1) stand-alone inverse optimization, (2) $K$-means clustering, (3) the $\MLIO_{\MIBP}$ model, (4) the benchmark model $\SEQMLIO$, and (5) the heuristics solution model $\EMBMLIO$ for the example case of clustering the given 80 decisions to three groups using the norm 2 distance between observed decisions and the optimal solution as the metric. \cref{fig:2D_Partition_comparison}(a) shows the feasible set and the observed decisions, and the other subfigures show the result from applying each model to the observed decisions and the given feasible set for finding three groups. \cref{fig:2D_Partition_comparison}(b) shows the clustering result of employing the conventional $K$-means algorithm to $\bX$. Note $K$-means is blind to $\Omega$. On the other hand, \cref{fig:2D_Partition_comparison}(c)  represents the outcome of applying the inverse learning model to the observed decisions while assuming (for comparison purposes) the observed decisions are homogeneous; the inverse learning model learns an optimal solution to $\FO$ that minimizes the defined metric for the case of only one cluster, because the inverse learning model, by definition, lacks partitioning capabilities.
The same effect as $K$-means is also seen in the benchmark model $\SEQMLIO$ in \cref{fig:2D_Partition_comparison}(d). Although the benchmark model $\SEQMLIO$ captures relatively acceptable optimal solutions (in the case of this representative example), the partitioning of the observed decisions suffers from being blind to $\Omega$. As shown in \cref{fig:2D_Partition_comparison}(e), using the modified heuristics method, $\EMBMLIO$, results in shifting the mappings found for different clusters based on the additional knowledge provided by $\Omega$. This finding is confirmed by the results from solving the $\MLIO_{\MIBP}$ model, which partitions the observed decisions in a manner that benefits from the given constraint data and is optimal based on the definition of $\MLIO$.  This feature is particularly important because $\EMBMLIO$ is solved in a fraction of the time required for solving $\MLIO_{\MIBP}$.  All these models---except $\MLIO_{\MIBP}$---solve the bi-dimensional example in a fraction of a second, including the heuristics models $\SEQMLIO$ and $\EMBMLIO$ with $\EMBMLIO$ finding the optimal solution that $\MLIO_{\MIBP}$ finds. As shown in the table, whereas $K$-means results in the least amount of total distance to the observed decisions, the centroids provided by $K$-means are not optimal solutions in general. 


\begin{table}[]
\begin{center}
\footnotesize
\renewcommand{\arraystretch}{1.2}
    \caption{Comparison of model performances for a representative bi-dimensional example in the case of three clusters. $\EMBMLIO$ finds the same solution that $\MIBP$ finds at optimality. $\EMBMLIO$ is solved within comparable times to the other solution models. $K$-means is not built for recovering the cost vectors, and both $K$-means and $\IO$ provide inferior optimality gap values based on the optimal solution found by $\MIBP$. Although the gap is slightly higher for $\SEQMLIO$, $\EMBMLIO$ reduces the gap to zero by finding the same optimal solution as $\MIBP$ in roughly 0.5 seconds.}
    \label{Table:2D_example_comparison}
    \begin{tabular}{>{\centering}p{0.12\textwidth}|>{\centering}p{0.1\textwidth}|>{\centering}p{0.1\textwidth}|>{\centering}p{0.18\textwidth}|>{\centering}p{0.2\textwidth}|>{\centering\arraybackslash}p{0.15\textwidth}}
    \hline\hline
\textbf{Model}        & \textbf{Learned Centroids} & \textbf{Recovered Cost Vectors} & \textbf{Total Optimality Gap}& \textbf{Total Distance to Observed Decisions}   & \textbf{Solution Time (Sec.)}     \\
\hline \hline
$K$-means & $\begin{matrix} (1.9,4.8)  \\ (5.4,7.3)  \\ (6.6,3.5)  \\ \end{matrix}$       &  NA    & 1.75 & 96.82  & 0.025      \\
\hline
$\IO$     & $\begin{matrix} \\ (6.7,7.3)  \\ \\ \end{matrix}$      &  $\begin{matrix}  (0.5,0.5) \end{matrix}$      & 10.86 & 315.1   & 0.121    \\
\hline \hline \hline
$\SEQMLIO$       & $\begin{matrix} (2,4.8)  \\ (5.4,7.6)  \\ (6.6,2)  \\ \end{matrix}$       & $\begin{matrix} (-1.0,0.0)  \\ (0.5,0.5)  \\ (0.0,-1.0)  \\ \end{matrix}$      & 0.32 & 136.7   & 0.188     \\
\hline
$\EMBMLIO$  & $\begin{matrix} (2,4.8)  \\ (5.5,7.4)  \\ (6.5,2)  \\ \end{matrix}$      & $\begin{matrix} (-1.0,0.0)  \\ (0.5,0.5)  \\ (0.0,-1.0)  \\ \end{matrix}$      & 0.00 & 135.2  & 0.594    \\
\hline
$\MIBP$ (optimal)    & $\begin{matrix} (2,4.8)  \\ (5.5,7.4)  \\ (6.5,2)  \\ \end{matrix}$    & $\begin{matrix} (-1.0,0.0)  \\ (0.5,0.5)  \\ (0.0,-1.0)  \\ \end{matrix}$    & 0.00 & 135.2 & $>$100  \\
    \hline  \hline
    \end{tabular}
\end{center}
\end{table}


Based on the results of the representative bi-dimensional example in this section, we observe that using the heuristic model $\EMBMLIO$ is significantly faster than and can be as capable as $\MLIO$ at finding optimal solutions for simple problems. In what follows, we apply the $\EMBMLIO$ model as a surrogate for $\MLIO$ in the setting of the diet recommendation problem for a set of non-homogeneous dietary behavior data. We compare the results of the $\EMBMLIO$ model with the benchmark model $\SEQMLIO$ and stand-alone unsupervised learning models such as $K$-means. We provide insights into the parallels between unsupervised clustering problems and the $\MLIO$ problem and compare the resulting recommended diets based on the known expert data.


\section{Data-Driven Diet Recommendation} \label{sec:application}

In this section, we apply the $\MLIO$ model and its two solution approaches to the NHANES dataset. In \cref{sec:app_data}, we briefly discuss the data used as the source for the observed dietary decisions and the constraint knowledge that forms the fixed and known $\Omega$ for the partially known optimization models. We then train models using $\EMBMLIO$, $\SEQMLIO$, and $K$-means methods to partition the observed decisions. In \cref{sec:app_training_results}, we provide comparisons on how these two models perform in partitioning the observed decisions to different groups and discuss relevant results. 

 The data-driven diet recommendation problem we consider contains different nutritional constraints that form the known and fixed $\Omega$ showcasing nutritional constraints that are extensively used in the literature \citep[see, e.g.,][]{garille2001stigler}. We base this problem on the nutritional constraints of a popular diet for reducing hypertension in adults, where the dietary preferences of the patients are unknown (cost vectors). We use a set of daily food intake data from patients to recover their dietary preferences and recover optimization problems such that each optimization problem is tailored to a group of similar patients based on their dietary behaviors. To do so, we employ the $\EMBMLIO$ model to partition the observed behaviors into a given number of clusters, recover the dietary preferences of the patients within each cluster, and recommend dietary choices for each group of similar patients. These recommendations will then have the desirable quality of minimizing in-group distances within each cluster and adhering optimally to the constraints of the DASH diet. For comparison, we consider the $K$-means clustering model and the $\SEQMLIO$ heuristics model as baseline benchmark models and base our performance metrics on the in-group distances of the observed decisions and behaviors of the patients in the same group with the recommended diets provided by the $\EMBMLIO$ and $\SEQMLIO$ models in the form of optimal solutions to optimization problems and their nutritional qualities. 

\subsection{Data} \label{sec:app_data}

We draw our analysis from several datasets. To incorporate the shared feasible set $\Omega$ that represents the known constraints in the diet recommendation problem, we use the recommendations of the Dietary Approaches to Stop Hypertension (DASH) eating plan \citep{dash_diet_2020} and form the constraints based on lower and upper bound nutritional bounds. These bounds are primarily based on evaluating the lower- and upper-bounds provided by the DASH eating plan for different food groups. Additionally, the forward optimization problem is capable of containing a myriad of other types of constraints (e.g., food-group-serving constraints), and the application of $\MLIO$ models is not contingent on the presence of specific types of constraints in $\Omega$. The observations for the application of the diet recommendation problem are individuals' daily food intakes, gathered from the open-access data from the National Health and Nutrition Examination Survey (NHANES)  dataset \citep{CDC_2020}. The reader is referred to the corresponding GitHub repository \citep{ahmadi2020opensource} for further details about this dataset.

Consider, from a more technical standpoint, the structure of $\FO$ in formulation \eqref{FO_i}, the variables of $\FO$ include the amount of daily intake for each food item, and the cost vector representing the preferences of the patients in each cluster. We allow non-homogeneity in observed decisions in terms of their dietary behavior and use the known constraints and the observed behaviors to cluster the decisions to homogeneous groups. The left-hand-side constraint matrix ($\bA$) is the amount of nutrients in a serving of each of the food items for each nutrient, and each nutrient has a lower- and upper-bound constraint. The right-hand-side constraint matrix ($\bb$) contains the bounds for each nutrient. These bounds are interpreted from the recommendation of the DASH eating plan \citep{dash_diet_2020}. 

For the analysis that follows, we consider a set of 900 patients with similar demographics and perform an 80-20 split of the data for training and testing purposes. Then, $\EMBMLIO$, $\SEQMLIO$, and $K$-means models are trained with 720 data points. We vary the number of clusters between 1 and 20 to observe the behavior of all models. We note the number of clusters indicates the number of different optimal solutions and cost vectors that the $\EMBMLIO$ and $\SEQMLIO$ models recover for the problem. We spotlight insights into how the $\EMBMLIO$ models approach the partitioning problem differently in comparison to $K$-means and discuss various features of the models.

\subsection{Results and Discussions} \label{sec:app_training_results}
In this section, we provide the results of training the $\EMBMLIO$ and the benchmark model $\SEQMLIO$ on the observed dietary behaviors of the patients described in the previous section. To compare these results against each other, we consider the results of the models for different food items (as the decision variables of the optimization problem) and the nutrient values (as the main constituents of $\Omega$) and compare them with the results of stand-alone applications of machine learning ($K$-means) clustering model ($K$-means). (Note that for the purposes of the analyses in this section, we use $\EMBMLIO$  instead of $\MLIO_{\MIBP}$ for the sake of computational efficiency.) 

\begin{figure}[h]
     \centering
     
    \begin{subfigure}{0.45\textwidth}
        \centering
        \includegraphics[width=\textwidth]{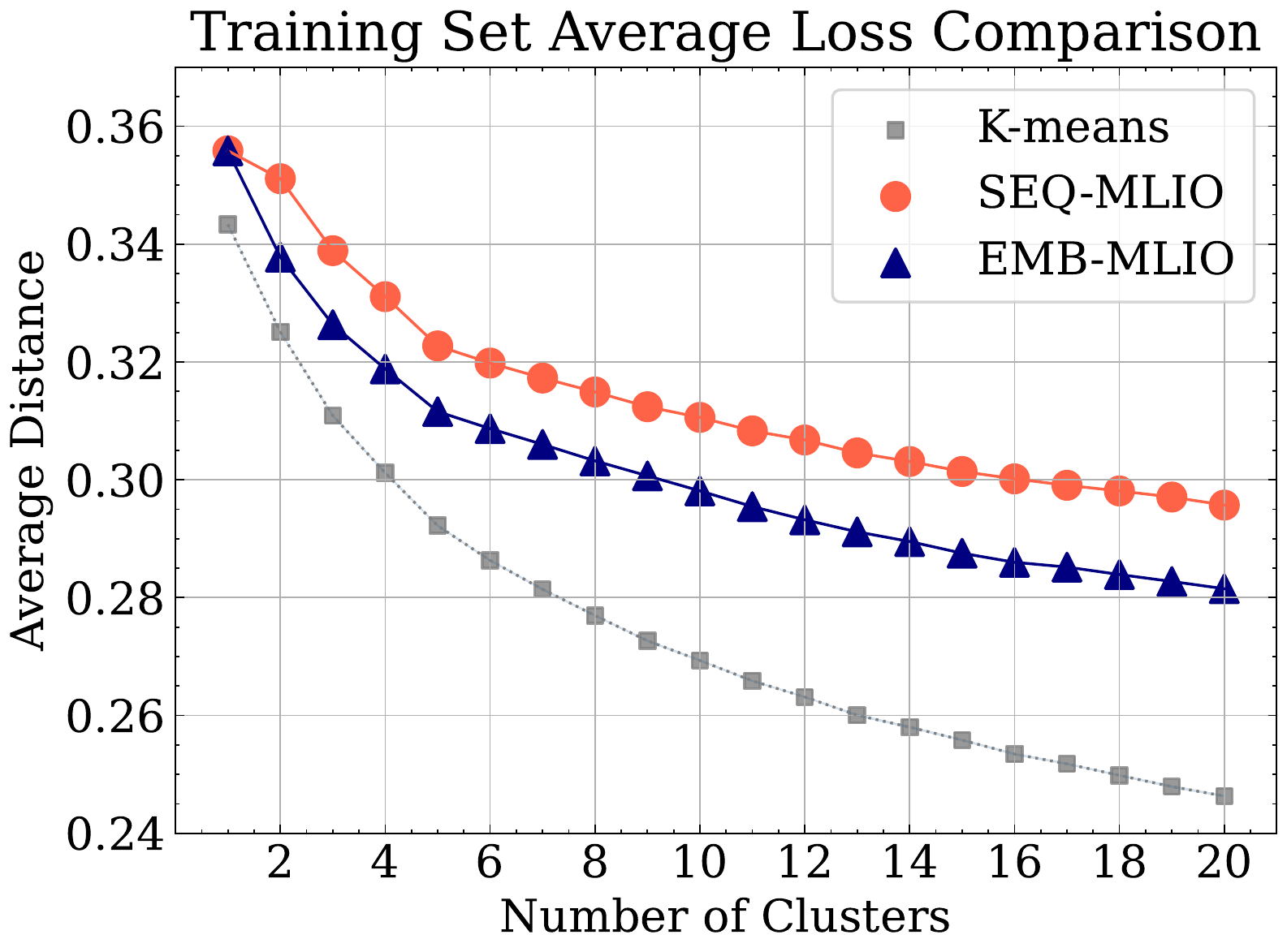}         \caption{}
        \label{fig:diet_training_comparison}
    \end{subfigure}
    \begin{subfigure}{0.45\textwidth}
        \includegraphics[width=\textwidth]{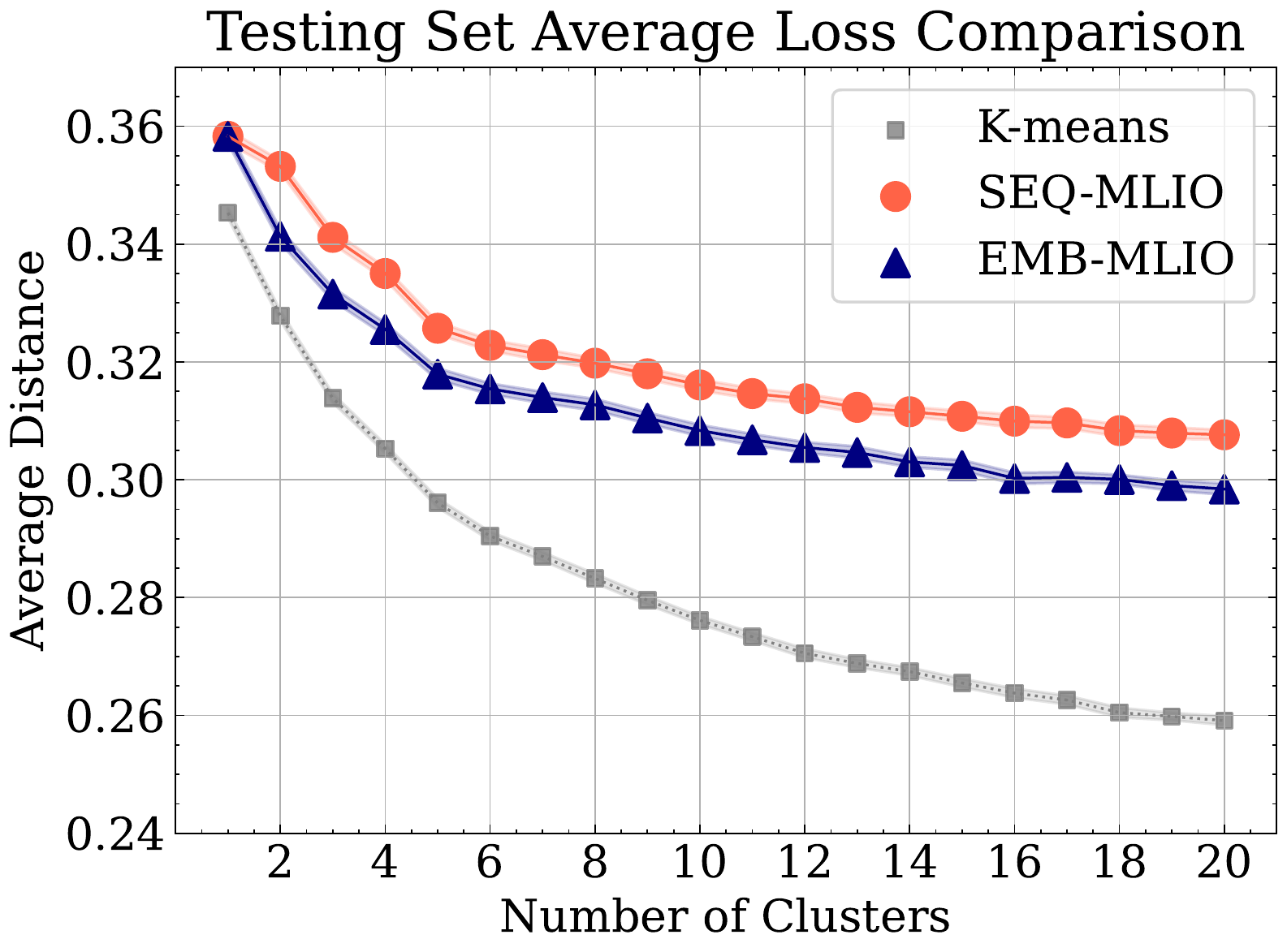}         \caption{} \footnotesize
        \label{fig:diet_testing_comparison}
    \end{subfigure}  
    \begin{subfigure}{0.45\textwidth}
        \includegraphics[width=\textwidth]{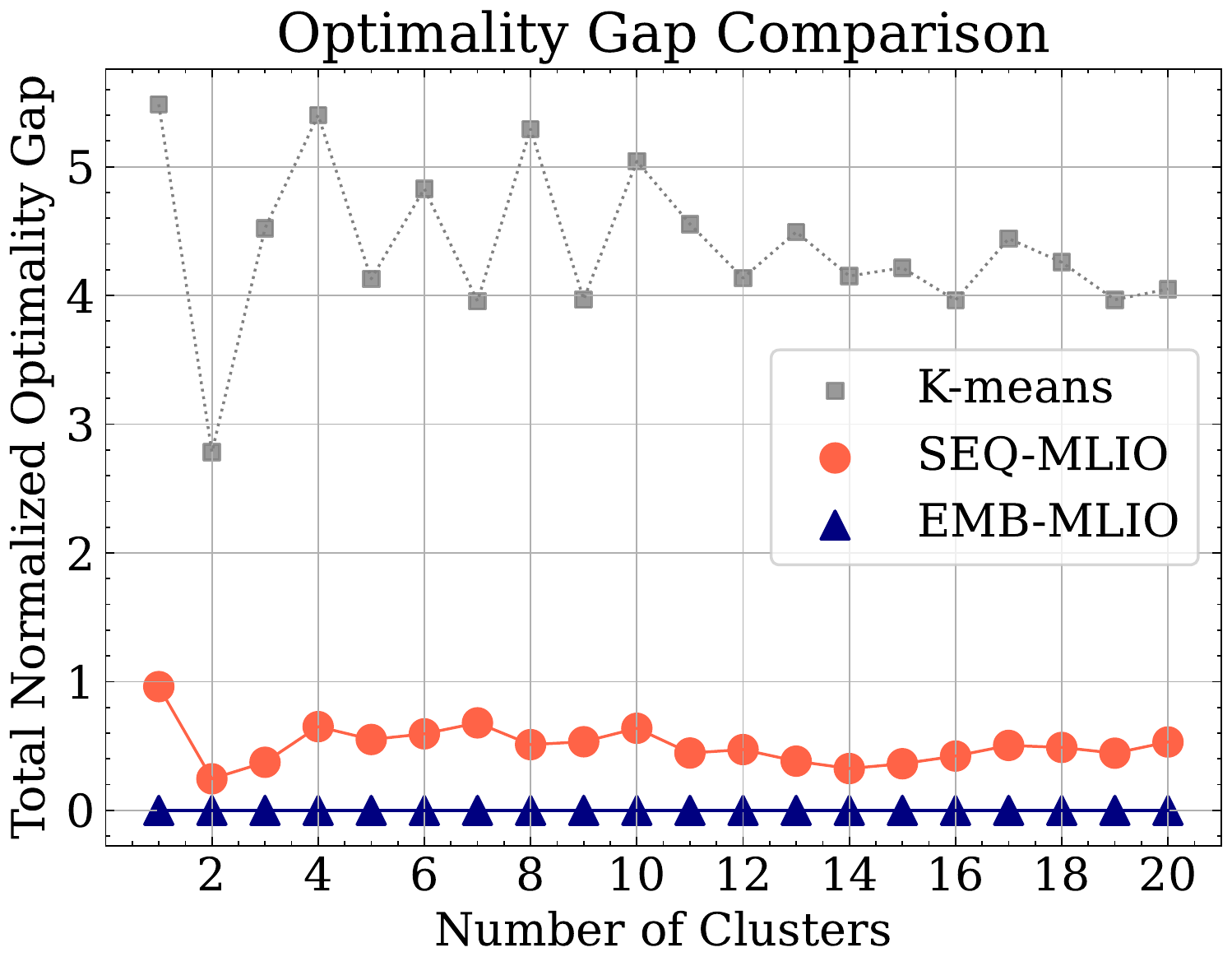}        \caption{}
        \label{fig:diet_optimality}
    \end{subfigure} 
    \hfill

    \caption{The average in-group distances are calculated for different numbers of clusters using the $K$-means (gray), $\SEQMLIO$ (orange), and $\EMBMLIO$ (navy). \cref{fig:diet_training_comparison} illustrates the average distance from the obtained solutions to the observed decisions within their clusters in the training data. \cref{fig:diet_testing_comparison} represents similar values with regard to the testing data. As shown, in both cases, $\EMBMLIO$ consistently achieves smaller distances than $\SEQMLIO$. \cref{fig:diet_optimality} shows the optimality gap of these models, where $\EMBMLIO$ achieves optimality and $\SEQMLIO$ averages around 0.5 in the total normalized optimality gap to the number of data points in the training set. In most cases, $K$-means represents a fluctuating pattern with normalized values above 4. For this analysis, training and testing samples include 720 and 180 data points of daily food intake, respectively.}\label{fig:diet_comparison}
\end{figure}

\cref{fig:diet_comparison} compares the  total in-group distances between observed decisions for the $\EMBMLIO$, $\SEQMLIO$, and $K$-means models. For different values of the numbers of clusters in the partition, $\EMBMLIO$ outperforms $\SEQMLIO$ for the training data. In comparison to $K$-means, whereas $\MLIO$ models perform on a subset of $\mathbb{R}^n$ and provide solutions that are farther from the observations than $K$-means, they provide optimal solutions over the known $\Omega$. \cref{fig:diet_comparison} shows the increase in the loss in comparison to $K$-means is not large. Additionally, the optimality-gap comparison figure shows how $K$-means does not provide optimal solutions with regard to $\Omega$ when proposing solutions for each cluster. We also note that, similar to conventional unsupervised learning models, this figure can also be used to get a notion of the appropriate number of clusters. In this case, because the curves are decreasing, one can either consider four clusters as an elbow point or consider 20 clusters for analysis. In what follows, we provide results for 20 clusters to facilitate comparisons across models.

\begin{figure}[]
\begin{center}
\includegraphics[width =1 \linewidth]{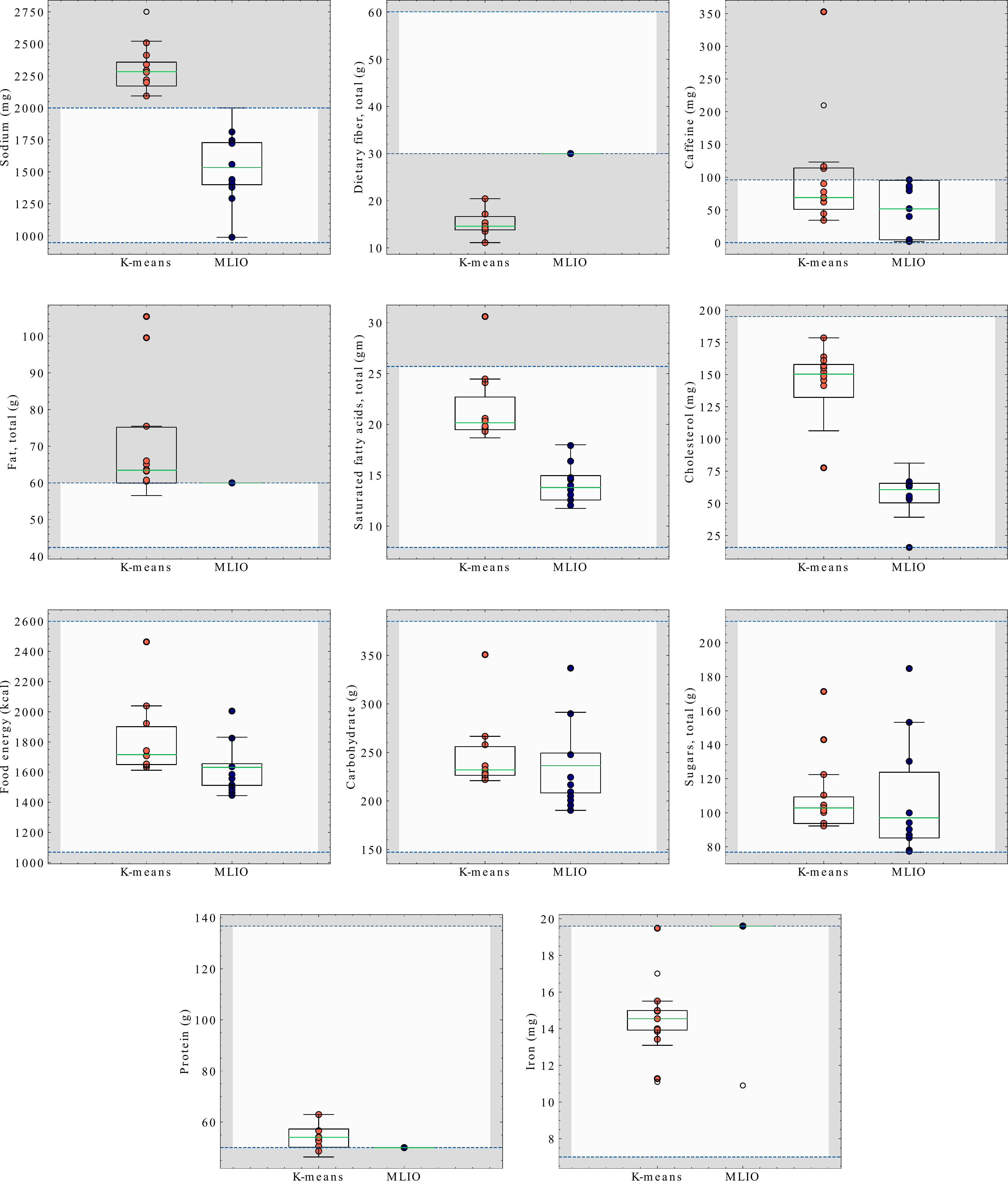}
\caption{\footnotesize   Box plots for the nutritional values of the diets recommended using the $K$-means and the $\EMBMLIO$ models. Nutrients with infeasible values are shaded. $K$-means recommends out of limit values for sodium and dietary fiber across the board, whereas $\EMBMLIO$ recommends diets that adhere to DASH limits (the non-shaded areas in the figures). $K$-means recommendations do not completely adhere to DASH limits for all clusters in the cases of total fat, protein, and cholesterol.} \label{Fig:Diet_nutrient_comparison}
\end{center}
\end{figure}

An important aspect of the diet recommendation problem is the nutritional quality of the diets generated by the algorithm. Thus, we  compare the results of the $\MLIO$ model with $K$-means models. As expected, because $K$-means models are naturally blind to the nutritional bounds, $K$-means models will only replicate the original behaviors of the patients, whereas $\MLIO$ will strike a balance between replicating behaviors and adhering to the hard constraints set by the DASH diet requirements. \cref{Fig:Diet_nutrient_comparison} compares the nutritional values of the recommended diet for both models for 20 clusters. For better representation, results are also shown via box plots representing the 10, 25, 50, 75, and 90 quantiles. Most notably, one can observe that MLIO restricts target nutrients of the DASH diet, such as sodium and cholesterol, but allows for comparable values with the original behaviors for other nutrients. The results show using $\MLIO$ allows dietitians to easily find diets that strike a balance between patient preferences and dietary goals. 

\begin{figure}[h]
     \centering
        \centering
        \includegraphics[width=1 \linewidth]{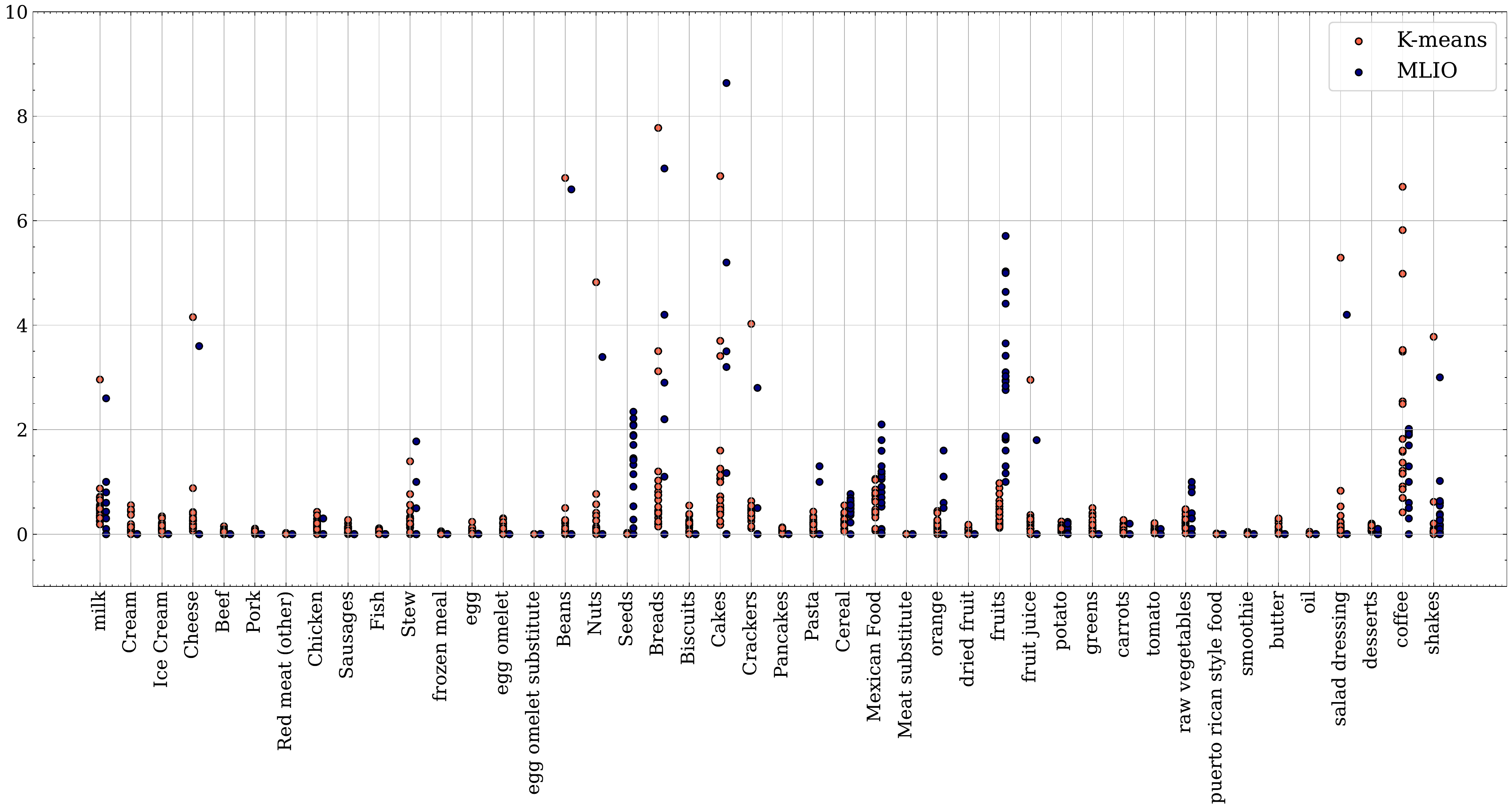}
        \caption{Comparison of  food group recommendations between the $\EMBMLIO$ and $K$-means models.}
        \label{fig:app_training_food_comparison}
\end{figure}

\begin{figure}[h]
     \centering
        \centering
        \includegraphics[width=1 \linewidth]{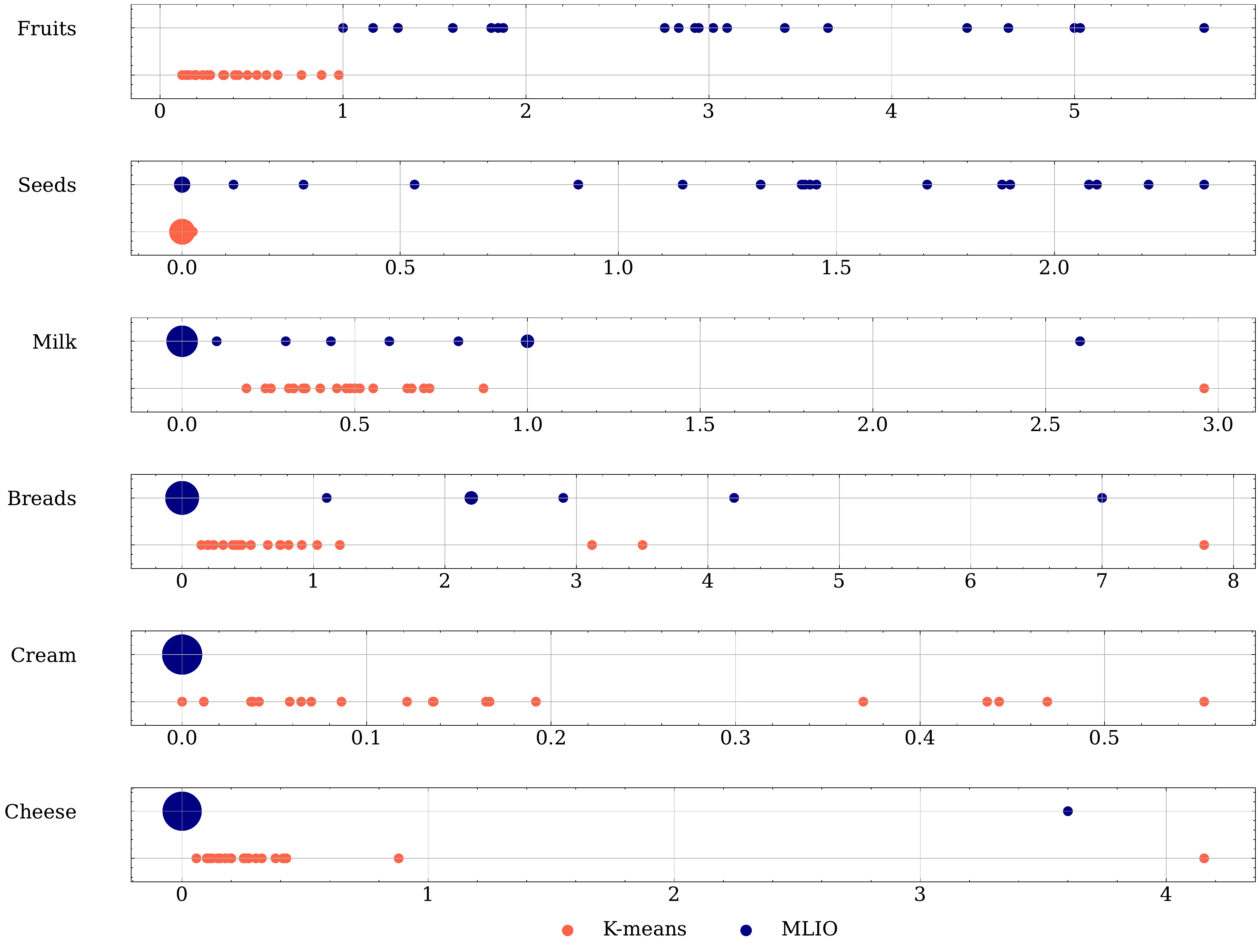}
        \caption{Individual food-group recommendations between the $\EMBMLIO$ and $K$-means models. $\EMBMLIO$ promotes healthy food groups such as fruits and seeds, whereas $K$-means replicates general unhealthy behaviors such as excessive consumption of creams in certain groups. Dots with increased radii indicate a higher number of groups that have that level of recommended intake for the food group.}
        \label{fig:app_training_food_comparison_select}
\end{figure}

In addition to comparing the in-group distances to the centroids and analyzing the quality of the recommended diets based on their nutritional values, we can analyze the centroids learned from the applications of the $\EMBMLIO$ model and a stand-alone unsupervised learning model such as  $K$-means. \cref{fig:app_training_food_comparison} compares the values of each food item for $\EMBMLIO$ (designated as $\MLIO$ in the figure) and $K$-means (designated as $K$-means in the figure). This comparison indicates the distinctions between the $K$-means and $\MLIO$ approaches in their learning solutions. Most notable among the food items, one can point out the fruit and seed food groups, which $\MLIO$ recommends significantly more than $K$-means. This difference in increased value represents adherence to nutritional bounds in recommending nutritious food groups. Additionally, as observed from \cref{fig:app_training_food_comparison_select}, recommendations for fruits are quite diverse for different clusters of patients. This can be attributed to the high bounds on dietary fiber. Because different clusters exhibit different levels of fiber, fruits are used primarily to adjust  diets to the required levels imposed by the DASH diet. On the other hand, items such as cheeses, creams, and dressings are among the foods that MLIO recommends noticeably less than $K$-means, which also represents restricting unhealthy groups, whereas  bread, crackers, and milk are among the foods that MLIO and $K$-means recommend to comparable extents, demonstrating MLIO's adherence to original behaviors as much as possible. The results of this figure confirm how $\MLIO$ balances between replicating original behaviors and incorporating constraint information.

The results provided in this section show the $\MLIO$ approach is capable of both (1) providing improved partitioning of observed decisions 
when information on healthy nutritional boundaries is available and also recovering optimization problems and (2) learning the best possible optimal solutions to such optimization problems based on the decisions partitioned in the same cluster. In the case of the diet recommendation, the recommendation system is capable of balancing between the healthy dietary constraints forming the dietary nutritional bounds and the observed behavior of the patients, and it provides diets that conform to the known dietary constraints. Additionally, we showed that, whereas $\MLIO$ restricts the unhealthy behaviors of patients, it is also capable of replicating their original preferences to the extent that the constraint-set information allows. 

\section{Concluding Remarks}\label{sec:conclusion}

In this paper, we present a novel data-driven approach to generating diet recommendations that incorporate both diverse dietary preferences and healthy dietary standards. This approach unifies inverse optimization and clustering. Under this hybrid approach, for the first time, inverse learning is embedded into a clustering scheme to recover unknown parameters and generate optimal solutions while ensuring optimality for generated solutions that act as cluster representatives. We demonstrated that sequentially applying unsupervised learning models (for partitioning of individuals) and inverse optimization models (for recovering the optimization problems) results in structurally distinct and suboptimal partitioning of observed decisions compared with solving the actual optimization problem; we proposed a solution method that provides some level of assurance over the partitioning quality.
To this end, we developed a novel solution strategy that \emph{embeds} clustering techniques and inverse optimization into each other. This strategy generalizes the clustering problem to the case where cluster representatives are subject to hard constraints and connects the generalized problem to existing inverse optimization models that recover the original optimization problems using observed data. 

The advantages of our hybrid approach are twofold. Results from applying the $\MLIO$ models to a bi-dimensional example with non-homogeneous observations over a fixed feasible set show how $\MLIO$ adapts to the known hard constraints of the problem setting. Equally important, applying $\MLIO$ models to the diet recommendation problem enables diet plan makers to incorporate group preferences by recovering missing objective functions and proposing optimal diets for different groups of individuals. Using $\MLIO$---as opposed to conventional machine learning techniques---allows dietary constraints to play a role in sharing how observations are clustered. By using $\MLIO$, we are able to redirect the recommendations of the systems toward diets that adhere to the nutritional limits set forth by experts, in our case, the DASH eating plans. 
Specifically for the diet recommendation problem,  a na\"ive clustering model results in the replication of unhealthy behaviors of patients.  By contrast, we show MLIO models are capable of recommending diets that comply with DASH nutritional standards. 



The MLIO approach described in this paper represents a step forward in modeling clustering problems with the goal of recommending optimal decisions for groups. This new approach will allow for data-driven decision support in applications such as the diet recommendation problem discussed in \cref{sec:application}. The MLIO approach can be used as the basis for developing an interactive recommendation tool that provides dietary recommendations based on their dietary history and dietary constraints catered to their needs.



%

%
%

\ACKNOWLEDGMENT{We gratefully acknowledge the financial support from the Johns Hopkins Discovery Award (2021--2023) and the Johns Hopkins Malone Center for Engineering in Healthcare Seed Grant (2020--2022). We appreciate the comments and suggestions from Fardin Ganjkhanloo, Julien Grand-Cl\'{e}ment, Huseyin Gurkan, Todd McNutt, and Nasrin Yousefi. We have also benefited from the feedback from seminar participants at the Johns Hopkins University's Center for Systems Science and Engineering and the Department of Applied Mathematics and Statistics, and session participants at 2020 INFORMS Annual Meeting, 2021 ACM Conference on Health, Inference, and Learning (CHIL), 2021 MSOM Annual Conference, and 2022 CORS/INFORMS International Conference.}


\clearpage

\iftrue
\begin{APPENDICES}{}
 
\section*{Appendix: Proofs}

\noindent  \emph{Proof of \cref{Prop:MLIO_feasibility}}
Let $\bX_1, \dots, \bX_L$ be some partition of $\bX_0$.  $\IO_1({\Omega,\bX_1})$, \dots, $\IO_L({\Omega,\bX_L})$ are all feasible and have at least one optimal solution set. Let $(\bc_1, \by_1, \bx^*_1, \bE_1), \hdots, (\bc_L, \by_L, \bx^*_L, \bE_L)$ be optimal solutions for $\IO_1({\Omega,\bX_1})$, \dots,$\IO_L({\Omega,\bX_L})$, respectively. Because $\bx^*_l$ are optimal for $\FO(\bc_l,\Omega)$ $\forall l \in \mathcal{L}$, the solution $(\bx^*_1, \dots, \bx^*_L$, $\bc_1, \dots, \bc_L )$ satisfies constraint \eqref{eq:opt_est_const}. Then,  $\bX_1, \dots, \bX_L$, $\bx^*_1, \dots, \bx^*_L$, $\bc_1, \dots, \bc_L$ is a feasible solution for $\MLIO({\bX_0, \Omega, L})$. \hfill \emph{Q.E.D.}
\smallskip

\noindent  \emph{Proof of \cref{Prop:MLIO_feasibility_all_partitions}}
Let $\bX_1, \dots, \bX_L \in \bP_L(\bX)$ be a partition of $\bX_0$ to $L$ clusters. Due to feasibility of $\IO$,  $\bx_1, \dots, \bx_L$ exist such that solution sets $(\left \{ \bc_1,\by_1,\bx_1,\bE_1 \right \}, \hdots, \left \{ \bc_L,\by_L,\bx_L,\bE_L \right \} $ are optimal  for $\IO (\bX_1,\Omega),\hdots,\IO (\bX_L,\Omega)$ respectively. As such, $\bX_1, \dots, \bX_L, \bx_1, \dots, \bx_L, \bc_1, \dots, \bc_L$ is feasible for $\MLIO({\bX_0, \Omega, L})$. Because the same is true for any partition in $\bP_L(\bX)$,  $\bP_L(\bX) \subseteq \hat{\Psi}_L$.
\hfill \emph{Q.E.D.}
\smallskip

\noindent  \emph{Proof of \cref{prop:MIBP_generalization_IL}}
Let $(\bx^*, \bv^*, \bY^*, \bC^*, \bE^*)$ be optimal for $\MIBP(\bX_0, \Omega, L)$. We first show that for each $l \in \mathcal{L}$, $(\bC^*_l, \bY^*_l, \bx^*_l, \hat{\bE^*_l})$ is feasible for $\IO (\bX^*_l,\Omega)$, where $(\bC^*_l, \bY^*_l, \bx^*_l, \hat{\bE^*_l})$ are defined similar to the statement of the proposition. This can be verified by noting constraints of $\MIBP$ are equivalent to $\IO$ for each value of $l$. Now, assume to the contrary that for some $l \in \mathcal{L}$, $(\bC^*_l, \bY^*_l, \bx^*_l, \hat{\bE^*_l})$ is not optimal for $\IO (\bX^*_l,\Omega)$. Then, some solution $(\bC'_l, \bY'_l, \bx'_l, \bE'_l)$ exist that is optimal for $\IO (\bX^*_l,\Omega)$, such that $\cD(\bX^*_l, \bx'_l) < \cD(\bX^*_l, \bx^*_l)$. However, we can build another solution to $\MIBP(\bX_0, \Omega, L)$ by replacing $(\bC^*_l, \bY^*_l, \bx^*_l, \hat{\bE^*_l})$ with  $(\bC'_l, \bY'_l, \bx'_l, \bE'_l)$ in the solution $(\bx^*, \bv^*, \bY^*, \bC^*, \bE^*)$. Because this new solution is feasible for $\MIBP(\bX_0, \Omega, L)$, this is a contradiction to $(\bx^*, \bv^*, \bY^*, \bC^*, \bE^*)$ being optimal for $\MIBP(\bX_0, \Omega, L)$ as we have found a new solution with a smaller objective value. \hfill \emph{Q.E.D.}
\smallskip

\noindent  \emph{Proof of \cref{Thm:MLIO_Equal}}
Let $(\bx^*, \bv^*, \bY^*, \bC^*, \bE^*)$ be optimal for $\MIBP(\bX_0, \Omega, L)$. Additionally, for $l \in \mathcal{L}$, let $\bX_l$ be the set of all elements of $\bX_0$ such that $\bx^k \in \bX_l$ if and only if $v^*_{k,l} = 1$. Considering that the solution  $(\bX^*, \bv^*, \bY^*, \bC^*, \bE^*)$ is feasible for $\MIBP(\bX_0, \Omega, L)$; it satisfies \eqref{MLIOclustering_constraint}. Therefore, we have $\bigcup_{l = 1}^{L} \bX_l  = \bX_0$, and as such, $\left \{ \bX_1,\hdots,\bX_L \right \}$ is a partition of $\bX_0$. Additionally, constraints \eqref{MLIOPrimalFeasblility1}, \eqref{MLIOStrongDual}, and \eqref{MLIODualFeas1} ensure for all $l \in \mathcal{L}$, $\bx^*_l$ is optimal for $\FO(\bC^*_l,\Omega)$. Therefore, $\forall l \in \mathcal{L}$, we have $\bx^*_l \in \Omega^{opt}(\bc_l)$ and as such, $(\left \{ \bX_1,\hdots,\bX_L \right \}, \bx^*, \bC^*)$ is feasible for $\MLIO$. Now, assume to the contrary that $(\left \{ \bX_1,\hdots,\bX_L \right \},\bx^*, \bC^*)$ is not optimal for $\MLIO$. Then, because any solution of $\MLIO$ can be mapped to at least one solution of $\MIBP$ by similar arguments as above, we would find a better solution for $\MIBP$ than $(\bx^*, \bv^*, \bY^*, \bC^*, \bE^*)$, which is a contradiction to $(\bx^*, \bv^*, \bY^*, \bC^*, \bE^*)$ being optimal for $\MIBP(\bX_0, \Omega, L)$. \hfill \emph{Q.E.D.}

\smallskip

\noindent  \emph{Proof of \cref{Prop:SEQ_MLIO_feasible_MLIO}}
Let $(\left \{ \bX_1,\hdots,\bX_L \right \}, \left \{ \bx^*_1,\hdots,\bx^*_L \right \}, \left \{ \bc_1,\hdots,\bc_L \right \})$ be the solution obtained from $\SEQMLIO$. To show $(\left \{ \bX_1,\hdots,\bX_L \right \},  \left \{ \bx^*_1,\hdots,\bx^*_L \right \}, \left \{ \bc_1,\hdots,\bc_L \right \})$ is feasible for $\MLIO$, we need to prove $\forall l \in \mathcal{L}$, $\bx^*_l \in \Omega^{opt}(\bc_l)$. However, notice that based on the definition of $\SEQMLIO$, $\forall l \in \mathcal{L}$, $\bx^*_l$ is part of an optimal solution to $\IO(\bX_l,\Omega)$, which is equivalent to $\bx^*_l$ being contained in $\Omega^{opt}(\bc_l)$. As such, the solution obtained by $\SEQMLIO$ is feasible for $\MLIO$.\hfill \emph{Q.E.D.}
\smallskip

\noindent  \emph{Proof of \cref{Prop:EMB_MLIO_feasible_MLIO}}
We first note that by definition of \cref{e_algorithm_for_MLIO}, for the solution set $(\left\{ \bX^i_1, \dots, \bX^i_L\right\}, \left\{\bx^i_1, \dots, \bx^i_L \right\}, \left \{ \bc^i_1,\hdots,\bc^i_L \right \})$ found in the $i^{th}$ iteration of \cref{e_algorithm_for_MLIO},  solutions $(\bc^i_1, \by^i_1, \bx^i_1, \bE^i_1), \hdots,(\bc^i_L, \by^i_L, \bx^i_L, \bE^i_L) $ exist that are optimal for $\IO({\Omega, \bX^i_1})$, ..., $\IO({\Omega, \bX^i_L})$, respectively. As such, recalling the definition of $\Omega^{opt}(\bc)$, we have $\bx^i_1 \in \Omega^{opt}(\bc^i_1), \hdots, \bx^i_L \in \Omega^{opt}(\bc^i_L)$. Because $\left\{ \bX^i_1, \dots, \bX^i_L\right\}$ is a partition of $\bX_0$, $(\left\{ \bX^i_1, \dots, \bX^i_L\right\}, \left\{\bx^i_1, \dots, \bx^i_L \right\}, \left \{ \bc^i_1,\hdots,\bc^i_L \right \})$ is feasible for $\MLIO$.\hfill \emph{Q.E.D.}
\smallskip

\noindent  \emph{Proof of \cref{Prop:EMB_MLIO_vs_SEQ_MLIO}}
Let $(\left\{ \bX^S_1, \dots, \bX^S_L\right\}, \left\{\bx^S_1, \dots, \bx^S_L \right\}, \left\{\bc^S_1, \dots, \bc^S_L \right\}, \left\{\bc^S_1, \dots, \bc^S_L \right\})$ be the $\MLIO$ solution resulting from the $\SEQMLIO$ approach. One can initialize \cref{e_algorithm_for_MLIO} using $\left\{ \bX^S_1, \dots, \bX^S_L\right\}$ as the initial partition of $\bX_0$. Considering that \cref{e_algorithm_for_MLIO} is strictly decreasing between iterations in terms of the total loss $\sum_{l=1}^{L} \mathcal{D} (\bx_l, \bX_l)$, if we denote the output partition and solution set of \cref{e_algorithm_for_MLIO} as $(\left\{ \bX^E_1, \dots, \bX^E_L\right\}, \left\{\bx^E_1, \dots, \bx^E_L \right\}, \left\{\bc^E_1, \dots, \bc^E_L \right\})$, we have $\sum_{l=1}^{L} \mathcal{D} (\bx^E_l, \bX^E_l) \leq \sum_{l=1}^{L} \mathcal{D} (\bx^S_l, \bX^S_l)$.\hfill \emph{Q.E.D.}
\smallskip

\noindent  \emph{Proof of \cref{Lemma_EMBMLIO}}
Assume to the contrary that $\exists$ $ \bar{\bX} = \left \{ \bX_1,\hdots,\bX_L \right \} \in \bP_L(\bX)$ that is visited twice in \cref{e_algorithm_for_MLIO} on iterations $i$ and $j$. For iterations $i$ and $j$, we have $\left\{ \bX^i_1, \dots, \bX^i_L\right\} = \left\{ \bX^j_1, \dots, \bX^j_L\right\}$. However, this means that for iterations $i$ and $j$, $\forall$ $l \in \mathcal{L}$, $\exists$ $(\bc^i_l, \by^i_l, \bx^i_l, \bE^i_l)$ optimal for $\IO({\Omega, \bX^i_l})$  and $(\bc^j_l, \by^j_l, \bx^j_l, \bE^j_l)$ optimal for $\IO({\Omega, \bX^j_l})$ where $\left\{\bx^i_1, \dots, \bx^i_L \right\} = \left\{\bx^j_1, \dots, \bx^j_L \right\}$. Therefore, for these iterations, we have $\sum_{l=1}^{L} \mathcal{D} (\bx^i_l, \bX^i_l) = \sum_{l=1}^{L} \mathcal{D} (\bx^j_l, \bX^j _l)$. But this is in contradiction to \cref{e_algorithm_for_MLIO} being strictly decreasing in terms of $\cD$ between iterations.  As such, each partition is visited at most once throughout \cref{e_algorithm_for_MLIO}.\hfill \emph{Q.E.D.}
\smallskip

\noindent \emph{Proof of \cref{prop:ALg_1_conv_iterations}}
Using the results of \cref{Lemma_EMBMLIO}, each potential partition $\left\{\bX_1, \dots, \bX_L\right\} \in \bP_L(\bX)$ is visited at most once in \cref{e_algorithm_for_MLIO}. Noting  the total number of possible partitions of $\bX_0$ to $L$ clusters is equal to $\left| \bP_L(\bX)\right|$, the number of iterations required for \cref{e_algorithm_for_MLIO} to stop and return a solution is at most $\left| \bP_L(\bX)\right|$. \hfill \emph{Q.E.D.}
\smallskip

\noindent  \emph{Proof of \cref{EMB_MLIO_convergence}}
Let $(\left\{\bX_1, \dots, \bX_L\right\}, \left\{\bx_1, \dots, \bx_L\right\},\left\{\bc_1, \dots, \bc_L\right\})$ be the final solution obtained by running \cref{e_algorithm_for_MLIO} on $\bX_0$. We first note that based on the operations done on lines 5 to 11 of \cref{e_algorithm_for_MLIO}, for each solution of $\EMBMLIO$ in each iteration, the first condition of \cref{Def:partial_optimal_solution} is satisfied as these lines find the best possible partition for a given set of optimal solutions in $\Omega^{opt}$. To investigate the second condition of \cref{Def:partial_optimal_solution}, we note that for $(\left\{\bX_1, \dots, \bX_L\right\}, \left\{\bx_1, \dots, \bx_L\right\},\left\{\bc_1, \dots, \bc_L\right\})$, three possible cases exist:
\begin{description}
    \item[Case 1.]  Stop criterion at line 15 of \cref{e_algorithm_for_MLIO}: In this case, solving $\IO$ for the last iteration solution did not yield a better optimal solution set. As such, the original optimal solution set $\left\{\bx_1, \dots, \bx_L \right\}$ was already the best possible solution set on $\Omega^{opt}$, which is equivalent to the second condition in \cref{Def:partial_optimal_solution}.
    \item[Case 2.] Stop criterion at line 19 of \cref{e_algorithm_for_MLIO}: In this case, although updating the partition did change $\bX_1, ... , \bX_L$ to a different partition $\bX^{New}_1, ... , \bX^{New}_L$, this has not changed the total loss and we have {$\bD (\bX^{New}_1,  \hdots , \bX^{New}_L, \bx_1,\hdots, \bx_L) = \bD (\bX_1,  \hdots , \bX_L, \bx_1,\hdots, \bx_L)$}. Therefore, although some other partition may generate the same total loss $\bD$ for the optimal solution set $\left\{\bx_1, \dots, \bx_L \right\}$,  no partition exists that provides a smaller total loss $\bD$. This is again equivalent to the second criterion of \cref{Def:partial_optimal_solution}.    
    \item[Case 3.]  Neither Stop criteria in \cref{e_algorithm_for_MLIO}: In this case, due to \cref{e_algorithm_for_MLIO} being strictly decreasing, the final solution found from \cref{e_algorithm_for_MLIO} is the global optimal solution because \cref{e_algorithm_for_MLIO} has exhausted all possible partitions in $\bP_L(\bX)$. Because the global optimal solution is a partial optimal solution as well, the statement holds. \hfill \emph{Q.E.D.}
\end{description} 
\smallskip

\end{APPENDICES}

\clearpage






\end{document}